\documentclass{amsart}
\usepackage{graphicx}
\usepackage{amsmath}
\usepackage{amssymb}
\usepackage{amsthm}

\theoremstyle{plain}
\newtheorem{theorem}{Theorem}[section]
\newtheorem{corollary}[theorem]{Corollary}

\newtheorem{lemma}[theorem]{Lemma}
\newtheorem{proposition}[theorem]{Proposition}
\newtheorem{definition}[theorem]{Definition}

\theoremstyle{definition}
\newtheorem{remark}[theorem]{Remark}

\numberwithin{equation}{section}

\begin{document}

\title{The weighted Yamabe problem with boundary}


\author{Pak Tung Ho}
\address{Department of Mathematics, Sogang University, Seoul, 04107, Korea}

\address{Korea Institute for Advanced Study, Seoul, 02455, Korea}

\email{paktungho@yahoo.com.hk}

\author{Jinwoo Shin}
\address{Korea Institute for Advanced Study, Hoegiro 85, Seoul 02455, Korea}
\email{shinjin@kias.re.kr}

\author{Zetian Yan}
\address{109 McAllister Building, Penn State University, University Park, PA 16802, USA}
\email{zxy5156@psu.edu}


\date{\today}

\begin{abstract}
We introduce a Yamabe-type flow
\begin{displaymath}
\left\{
\begin{array}{ll}
\frac{\partial g}{\partial t} &=(r^m_{\phi}-R^m_{\phi})g \\
\frac{\partial \phi}{\partial t} &=\frac{m}{2}(R^m_{\phi}-r^m_{\phi})
\end{array}
\right. ~~\mbox{ in }M ~~\mbox{ and }~~
H^m_{\phi}=0 ~~\mbox{ on }\partial M
\end{displaymath}
on a smooth metric measure space with boundary $(M,g, v^mdV_g,v^mdA_g,m)$, where $R^m_{\phi}$ is the associated weighted scalar curvature, $r^m_{\phi}$ is the average of the weighted scalar curvature, and
$H^m_{\phi}$ is the weighted mean curvature. We prove the long-time existence and convergence of this flow.
\end{abstract}

\maketitle

\section{Introduction}\sloppy

Suppose $M$ is a compact, $n$-dimensional manifold
without boundary, where $n\geqslant 3$, and $g_0$ is a Riemannian metric on $M$.
As a generalization of the Uniformization Theorem, the Yamabe problem is to find a
metric conformal to $g_0$ such that its scalar curvature is constant.
This was first introduced by Yamabe \cite{Yamabe},
and was solved by Aubin \cite{Aubin0}, Trudinger \cite{Trudinger} and Schoen \cite{Schoen}.

The Yamabe flow is defined as
$$
\frac{\partial g}{\partial t}=(r_g-R_{g})g,
$$
where $R_g$ is the scalar curvature of $g$ and $r_g$ is the average of $R_g$:
$$r_g=\frac{\int_MR_gdV_g}{\int_MdV_g}.$$
This was first introduced by  Hamilton in \cite{Hamilton}.
Hamilton conjectured that, for every initial metric, the flow converges to a conformal metric of constant scalar curvature. In the case when $Y(M, g_0)\leqslant 0$, it is not difficult to show that the conformal factor is uniformly bounded above and below. Moreover, the flow converges to a metric of constant scalar curvature as $t\to \infty$.

The case $Y(M, g_0)>0$ is more interesting. Chow \cite{Chow} proved the convergence of the flow for locally conformally flat metrics with positive Ricci curvature. Ye \cite{Ye} later extended the result to all locally conformal flat metrics. Later, Brendle \cite{Brendle4} proved convergence of the flow for all conformal classes and arbitrary initial metrics when $3\leqslant n\leqslant 5$, and extended the results to higher dimensions \cite{Brendle5}.

Now consider a compact, $n$-dimensional manifold $M$
with smooth boundary $\partial M$, where $n\geqslant 3$,
and $g_0$ is a Riemannian metric on $M$.
One can still talk about the Yamabe problem
for manifold with boundary, and there are two types.
For the first type, one would like to find a conformal metric $g$
such that its scalar curvature $R_g$ is constant in $M$ and
its mean curvature $H_g$ is zero on $\partial M$.
For the second type, one would like to find a
conformal metric $g$
such that its scalar curvature $R_g$ is zero in $M$ and
its mean curvature $H_g$ is constant on $\partial M$.
These problems have studied by many authors.
See \cite{Escobar2,Escobar1} for example.

Similar to the Yamabe flow,
Brendle  introduced some geometric flows in \cite{Brendle}
to study the Yamabe problem for manifolds with boundary.
For the first type, the geometric flow is defined as
\begin{equation}\label{3}
\frac{\partial g}{\partial t}=-(R_g-r_{g})g\mbox{ in }M~~\mbox{ and }~~H_{g}=0\mbox{ on }\partial M.
\end{equation}
Almaraz and Sun has considered in \cite{Almaraz&Sun} the convergence of the flow (\ref{3}). On the other hand, for the second type, the geometric flow is defined as
\begin{equation}\label{4}
\frac{\partial g}{\partial t}=-(H_{g}-h_{g})g\mbox{ on }\partial M~~\mbox{ and }~~R_{g}=0\mbox{ in }M,
\end{equation}
where $h_{g}$ is the average of the mean curvature $H_{g}$:
\begin{equation*}
h_{g}=\frac{\int_{\partial M} H_{g}dA_{g}}{\int_{\partial M}dA_{g}}.
\end{equation*}
In \cite{Almaraz}, Almaraz has studied the convergence of the flow (\ref{4}). See also
\cite{Chen&Ho,Ho&Lee&Shin,Ho&Shin} for results related to the flows (\ref{3}) and (\ref{4}).

In this paper, in the same spirit of \cite{Yan}, we generalize the Yamabe flow to smooth metric measure spaces with boundary. To explain the results of this paper,
we require some terminology.

\begin{definition}
Let $(M,\partial M, g)$ be a Riemannian manifold with boundary $\partial M$ and let us denote by $dV_g$ and $dA_g$ the volume element induced by $g$ on $M$ and $\partial M$, respectively. Set a function $\phi\in C^{\infty}(M)$ and a dimensional parameter $m\in [0,\infty)$. In the case $m=0$, we require that $\phi=0$.
A smooth metric measure space with boundary is
a five-tuple $(M,g,e^{-\phi}dV_g,e^{-\phi}dA_g,m)$. We frequently denote a smooth metric measure space by $(M,g, v^mdV_g,v^mdA_g,m)$ where $\phi$ and $v$ are related by
$e^{-\phi}=v^m$.
\end{definition}
The \textit{weighted scalar curvature} $R^m_{\phi}$ of a smooth metric measure space with boundary $(M,g,e^{-\phi}dV_g,e^{-\phi}dA_g,m)$ is
\begin{equation}\label{defn1}
R^m_{\phi}:=R_g+2\Delta_g \phi-\frac{m+1}{m}|\nabla\phi|_g^2,
\end{equation}
where $R_g$ and $\Delta_g$ are the scalar curvature and the Laplacian associated to the metric $g$, respectively. The \textit{weighted mean curvature} is
\begin{equation}
H^m_{\phi}=H_{g}+\frac{\partial\phi}{\partial\nu_{g}},
\end{equation}
where $H_{g}$ and $\displaystyle\frac{\partial}{\partial\nu_{g}}$
are the mean curvature and the outward normal derivative with respect to $g$, respectively.

Conformal equivalence between smooth metric measure spaces are defined as follows, see \cite{Case15} for more details.
\begin{definition}\label{condef}
Smooth metric measure spaces with boundary $(M, g, e^{-\phi}dV_g,e^{-\phi}dA_g, m)$ and $(M, \hat{g}, e^{-\hat{\phi}}dV_{\hat{g}},e^{-\hat{\phi}}dA_{\hat{g}}, m)$ are pointwise conformally equivalent if there is a function $\sigma\in C^{\infty}(M)$ such that
\begin{equation}\label{conformal}
(M, \hat{g}, e^{-\hat{\phi}}dV_{\hat{g}},e^{-\hat{\phi}}dA_{\hat{g}}, m)=(M,e^{\frac{2\sigma}{m+n-2}}g, e^{\frac{(m+n)\sigma}{m+n-2}}e^{-\phi}dV_g,e^{\frac{(m+n)\sigma}{m+n-2}}e^{-\phi}dA_g, m).
\end{equation}
In the case $m=0$, conformal equivalence is defined in the classical sense.
\end{definition}
If we denote $e^{\frac{1}{2}\sigma}$ by $w$, (\ref{conformal}) is equivalent to
\begin{equation}\label{1.5}
(M, \hat{g}, e^{-\hat{\phi}}{\rm dvol}_{\hat{g}},e^{-\hat{\phi}}dA_{\hat{g}}, m)=(M,w^{\frac{4}{m+n-2}}g,w^{\frac{2(m+n)}{m+n-2}}e^{-\phi}dV_{g},w^{\frac{2(m+n)}{m+n-2}}e^{-\phi}dA_{g}, m),
\end{equation}
which is an alternative way to formulate the conformal equivalence of smooth metric measure spaces.

\begin{definition}
Let $(M, g, e^{-\phi}dV_g, m)$ be a smooth metric measure space. The weighted Laplacian $\Delta_{\phi}: C^{\infty}(M)\to C^{\infty}(M)$ is the operator defined as
\begin{displaymath}
\Delta_{\phi}\psi=\Delta\psi-\langle\nabla \phi, \nabla \psi\rangle_g ~~\mbox{ for any }\psi\in C^\infty(M),
\end{displaymath}
It is formally self-adjoint with respect to the measure $e^{-\phi}dV_{g}$. For more about smooth metric measure spaces, we refer the readers to \cite{Case15,Case12,Case13,Han}.
\end{definition}

\begin{definition}
Given a smooth metric measure spaces with boundary
$(M, g, e^{-\phi}dV_g,e^{-\phi}dA_g, m)$, the weighted conformal Laplacian $(L_{\phi}^m, B_{\phi}^m)$ is given by the interior operator and boundary operator
\begin{equation}\label{2.3}
\begin{split}
L_{\phi}^m&=-\Delta_{\phi}+\frac{n+m-2}{4(n+m-1)}R^m_{\phi} ~~\mbox{ in }M,\\
B_{\phi}^m&=\frac{\partial}{\partial\nu_{g}}+\frac{n+m-2}{2(n+m-1)}H^m_{\phi} ~~\mbox{ on }\partial M,
\end{split}
\end{equation}
where $\nu_g$ is the outward unit normal with respect to $g$.
\end{definition}

Consequently, in the formulation of (\ref{1.5}), the transformation law of the weighted scalar curvature and the weighted mean curvature \cite[Proposition 1]{Posso18} are
\begin{equation}\label{1.6}
\begin{split}
R^m_{\phi}&=\frac{4(n+m-1)}{n+m-2}w^{-\frac{m+n+2}{m+n-2}}L^m_{\phi_0}w ~~\mbox{ in }M,\\
H^m_{\phi}&=\frac{2(n+m-1)}{n+m-2}w^{-\frac{n+m}{n+m-2}}B_{\phi_0}^mw~~\mbox{ on }\partial M.
\end{split}
\end{equation}

Given a compact smooth metric measure space without boundary $(M,g,e^{-\phi}dV_{g}, m)$,
the \textit{weighted Yamabe problem} is to find
another smooth metric measure space $(M,\hat{g},e^{-\hat{\phi}}dV_{\hat{g}}, m)$
conformal to $(M,g,e^{-\phi}dV_{g}, m)$ such that
its weighted scalar curvature $R^m_{\hat{\phi}}$ is constant.
This was first introduced and studied by Case in \cite{Case15}.
See also \cite{Case12,Case19,Souza20,Posso21} for results related to the weighted Yamabe problem.

Similarly, the weighted Yamabe problem with boundary is to find $(M, \hat{g}, e^{-\hat{\phi}}dV_{\hat{g}},e^{-\hat{\phi}}dA_{\hat{g}}, m)$ conformal to $(M, g, e^{-\phi}dV_g,e^{-\phi}dA_g, m)$ such that $R_{\hat{\phi}}^m$ is constant in $M$
and $H_{\hat{\phi}}^m$ is zero on $\partial M$. In view of (\ref{1.6}),
it is equivalent to solve
\begin{equation}\label{2.4}
\begin{split}
L_{\phi}^mw&=\frac{n+m-2}{4(n+m-1)}\lambda w^{\frac{n+m+2}{n+m-2}}~~\mbox{ in }M,\\
B_{\phi}^mw&=0~~\mbox{ on }\partial M
\end{split}
\end{equation}
for some constant $\lambda$. This has been studied by Posso in \cite{Posso18}.

In the spirit of \cite{Brendle}, we introduce a Yamabe-type flow on the smooth metric measure space with boundary
$(M, g, e^{-\phi}dV_g,e^{-\phi}dA_g, m)$, $m\in (0,\infty)$, which is a natural way to solve the weighted Yamabe problem with boundary. The definition of the flow arises from the following observation in \cite{Yan}. In the sense of Definition \ref{condef}, the metric $(e^{\phi})^{\frac{2}{m}}g$ is fixed within the conformal class of $(M, g, e^{-\phi}dV_g,e^{-\phi}dA_g, m)$. Based on this observation, we define the (normalized) \textit{weighed Yamabe flow with boundary} as
\begin{equation}\label{flow}
\left\{
\begin{array}{ll}
\frac{\partial g}{\partial t} &=(r^m_{\phi}-R^m_{\phi})g \\
\frac{\partial \phi}{\partial t} &=\frac{m}{2}(R^m_{\phi}-r^m_{\phi})
\end{array}
\right. ~~\mbox{ in }M ~~\mbox{ and }~~
H^m_{\phi}=0 ~~\mbox{ on }\partial M,
\end{equation}
where $r_{\phi}^m$ is the average of the weighted scalar curvature $R_{\phi}^m$; i.e.
\begin{equation}\label{3.3}
r_{\phi}^m=\frac{\int_M R_{\phi}^m e^{-\phi}dV_{g}}
{\int_M  e^{-\phi}dV_{g}}.
\end{equation}
On the one hand, equation (\ref{flow}) is analogous to the Yamabe flow with boundary (\ref{3}). 
On the other hand,  the flow (\ref{flow}) is ``sub-critical"
in the sense that $\frac{2(n+m)}{n+m-2}<\frac{2n}{n-2}$. As a result, we can establish the sequential compactness in Proposition \ref{prop5.1}, which is 
the main difference between our flow (\ref{flow}) and the geometric flow (\ref{3}) (see \cite{Almaraz&Sun} for more details). Moreover, we adapt an argument of Brendle \cite{Brendle4} to establish long-time existence and convergence of the weighted Yamabe flow with boundary.

\begin{theorem}\label{main}
On a smooth metric measure space with boundary $(M, g, e^{-\phi}dV_g,e^{-\phi}dA_g, m)$, where $(M,\partial M, g)$ is a compact Riemannian manifold with boundary of dimension $n\geqslant3$, for every choice of the initial metric and the measure, the weighted Yamabe flow with boundary (\ref{flow}) exists for all time and converges to a metric with constant weighted scalar curvature.
\end{theorem}

This paper is organized as follows. As mentioned above, we first deal with the positive case.

In Section 2, we prove that the conformal factor $w(t)$ cannot blow up in finite time by bounding $w(t)$ from above and below on any finite time interval $[0,T]$. The long-time existence follows from this.

Convergence of the flow (\ref{flow}) will be based on the crucial observation in Proposition \ref{prop4.3}. In Section 3,
by assuming Proposition \ref{prop4.3}, we   obtain decay rates of $r^m_{\phi(t)}$ and the uniform upper bound of $|R^m_{\phi(t)}-r^m_{\phi(t)}|$ in $L^2$ norm. Together with the interior regularity theorem and estimates on the boundary, we show that $w(t)$ is uniformly bounded above and below on $[0,\infty)$, which implies
the convergence of the weighted Yamabe flow with boundary.

In Section 4, we complete the proof of Proposition \ref{prop4.3} by using the spectral theorem of self-adjoint operators and asymptotic analysis.

In Section 5, in the same spirit as \cite{Ye}, we refine the argument in Section 2 to obtain the uniform bound on $w(t)$ and prove the long-time existence and smooth convergence in the negative case. Besides, in the zero case, we obtain the Harnack inequality such that uniform smooth estimates hold.

\section{Long time existence}

In this section we collect some basic facts for smooth metric measure spaces and prove various properties of the weighted Yamabe flow with boundary that will be used throughout this paper.

\begin{definition}\label{2.1}
On a smooth metric measure space with boundary
$(M, g, e^{-\phi}dV_g,e^{-\phi}dA_g, m)$ which is conformal to $(M, g_0, e^{-\phi_0}dV_{g_0},e^{-\phi_0}dA_{g_0}, m)$ in the formulation of (\ref{1.5}), analogous to the classical Yamabe problem with boundary, we define the normalized energy  $E(w)$ as
\begin{equation*}
E_{g_0,\phi_0}(w)=\frac{\int_MwL_{\phi_0}^m(w)e^{-\phi_0}dV_{g_0}+\int_MwB_{\phi_0}^m(w)e^{-\phi_0}dA_{g_0}}{(\int_Mw^{\frac{2(n+m)}{n+m-2}}e^{-\phi_0}dV_{g_0})^{\frac{n+m-2}{n+m}}}.
\end{equation*}
\end{definition}
\begin{remark} By the transformation law in (\ref{1.6}), the normalized energy $E_{g_0,\phi_0}(w)$ can be written as
\begin{displaymath}
E_{g_0,\phi_0}(w)=\frac{n+m-2}{4(n+m-1)}\frac{\int_MR_\phi^me^{-\phi}dV_g+2\int_{\partial M}H_\phi^me^{-\phi}dA_g}{(\int_M e^{-\phi}dV_g)^{\frac{n+m-2}{n+m}}}.
\end{displaymath}
\end{remark}
We set
\begin{equation*}
Y_{n,m}[(g_0, \phi_0)]=\inf\left\{E_{g_0,\phi_0}(w): 0<w\in C^\infty(M)\right\}.
\end{equation*}
In order to analyze the long time behavior of the solutions of (\ref{flow}), we consider three different cases:
\begin{equation*}
    \begin{array}{ll}
    ~~\mbox{Positive case}~~: & Y_{n,m}[(g_0, \phi_0)]>0,  \\
    ~~\mbox{Zero case}~~: & Y_{n,m}[(g_0, \phi_0)]=0,  \\
    ~~\mbox{Negative case}~~: & Y_{n,m}[(g_0, \phi_0)]<0.
\end{array}
\end{equation*}
Similar to \cite[Lemma 2.1]{Brendle}, we have the following lemma. 
\begin{lemma}\label{threecases}
There exists a smooth metric measure space with boundary
$(M, g, e^{-\phi}dV_g,e^{-\phi}dA_g, m)$ which is conformal to $(M, g_0, e^{-\phi_0}dV_{g_0},e^{-\phi_0}dA_{g_0}, m)$ such that
\begin{equation*}
    R^m_{\phi}>0~(=0, <0\mbox{ respectively}) ~~\mbox{ in }M ~~\mbox{ and }~~
H^m_{\phi}=0 ~~\mbox{ on }\partial M,
\end{equation*}
if $Y_{n,m}[(g_0, \phi_0)]>0 ~(=0, <0\mbox{ respectively})$.
\end{lemma}

In light of the discussion in \cite{Ye}, in the case $Y_{n,m}[(g_0, \phi_0)] \leqslant 0$, it is not difficult to show convergence of the flow (\ref{flow}) as $t\to \infty$.  We postpone the proof to Section 5 and deal with the positive case first.

Hereafter, we choose $(M, g_0, e^{-\phi_0}dV_{g_0},e^{-\phi_0}dA_{g_0}, m)$ to be the initial metric measure space with $Y_{n,m}[(g_0, \phi_0)]> 0$. 
By conformal change, we may assume that
the initial weighted mean curvature $H_{\phi_0}^m$ satisfies
\begin{equation}\label{3.1}
H^m_{\phi_0}=0\mbox{ on }\partial M.
\end{equation}
Since the weighted Yamabe flow preserves the conformal structure, we may write \begin{equation}
\begin{cases}
&g(t)=w(t)^{\frac{4}{n+m-2}}g_0, \\
&e^{-\phi(t)}=w(t)^{\frac{2(m+n)}{n+m-2}}e^{-\phi_0},
\end{cases}
\end{equation}
as the solution of (\ref{flow}) with $(g(0), \phi(0))=(g_0, \phi_0)$. Hence, the first equation of (\ref{flow}) reduces to the following evolution equation for the conformal factor
\begin{equation}\label{3.5}
\frac{\partial}{\partial t}w(t)^{\frac{n+m+2}{n+m-2}}=\frac{m+n-2}{4}\left(\frac{4(n+m-1)}{n+m-2}\Delta_{\phi_0}w-R_{\phi_0}^mw+r^m_{\phi(t)}w^{\frac{n+m+2}{n+m-2}}\right)~~\mbox{ in }M.
\end{equation}
It follows from (\ref{3.1}) and (\ref{2.3}) that
the second equation in (\ref{flow}) is equivalent to
\begin{equation}\label{3.6}
\frac{\partial w(t)}{\partial\nu_{g_0}}=0~~\mbox{ on }\partial M.
\end{equation}
Hence the conformal factor $w(t)$ satisfies the evolution equations
\begin{equation}\label{wevolution}
    \left\{
\begin{array}{ll}
    \frac{\partial w(t)}{\partial t}&=-\frac{m+n-2}{4}(R^m_{\phi(t)}-r^m_{\phi(t)})w(t)~~\mbox{ in }M,\\
    \frac{\partial w(t)}{\partial\nu_{g_0}}&=0~~\mbox{ on }\partial M.
\end{array}
\right.
\end{equation}

By direct calculation, integration by parts on a smooth metric measure space with boundary $(M, g, e^{-\phi}dV_g,e^{-\phi}dA_g, m)$ takes the following form
\begin{equation}\label{parts}
    \begin{split}
      \int_M \langle \nabla f,X\rangle e^{-\phi}dV_g=&-\int_M f {\rm div}_\phi(X)e^{-\phi}dV_g
      +\int_{\partial M} f\langle X,\nu_g\rangle e^{-\phi}dA_g
    \end{split}
\end{equation}
for any smooth vector field $X$ in $M$, where ${\rm div}_\phi(X)={\rm div} X-\langle X, \nabla \phi\rangle$.

Since
\begin{equation}\label{vol}
    \frac{d}{d t}\int_M e^{-\phi(t)}dV_{g(t)}
 =\frac{n+m}{2}\int_M (r^m_{\phi}-R^m_{\phi})e^{-\phi}dV_{g}=0,
\end{equation}
we may assume that
\begin{equation}\label{normalize1}
\int_M e^{-\phi(t)}dV_{g(t)}=1
\end{equation}
for all $t\geqslant 0$. With this normalization, the average of the weighted scalar curvature can be written as
\begin{equation}\label{normalize2}
r^m_{\phi}(t)=\int_M R^m_{\phi(t)} e^{-\phi(t)}dV_{g(t)}.
\end{equation}

By (\ref{wevolution}), differentiating (\ref{3.6}) with respect to $t$ yields $\displaystyle
\frac{\partial R_{\phi(t)}^m}{\partial\nu_{g_0}}=0$ on $\partial M$. Since $\displaystyle\frac{\partial}{\partial\nu_{g(t)}}=w(t)^{-\frac{2}{n+m-2}}\frac{\partial}{\partial\nu_{g_0}}$,
this is equivalent to
\begin{equation}\label{3.8}
\frac{\partial R_{\phi(t)}^m}{\partial\nu_{g(t)}}=0~~\mbox{ on }\partial M.
\end{equation}
Combining with (1.5) in \cite{Yan}, we deduce that the weighted scalar curvature satisfies the following evolution equations
\begin{equation}\label{Revolution}
    \left\{
     \begin{array}{ll}
       \frac{\partial R_{\phi(t)}^m}{\partial t}&=(n+m-1)\Delta_{\phi(t)}R_{\phi(t)}^m+R_{\phi(t)}^m(R_{\phi(t)}^m-r_{\phi(t)}^m)~~\mbox{ in }M,\\
       \frac{\partial R_{\phi(t)}^m}{\partial\nu_{g(t)}}&=0~~\mbox{ on }\partial M.
      \end{array}
    \right.
\end{equation}

Using the evolution equation (\ref{Revolution}), we obtain
\begin{equation}\label{eq:r}
\frac{d}{d t}r^m_{\phi(t)}
=-\frac{n+m-2}{2}\int_M (r^m_{\phi(t)}-R^m_{\phi(t)})^2e^{-\phi(t)}dV_{g(t)}\leqslant 0.
\end{equation}
Observe that $r^m_{\phi}(t)>0$ since $Y_{n,m}[(g, \phi)]>0$. Hence, $r^m_{\phi(t)}$ is bounded above and below, i.e.
\begin{equation}\label{boundr}
    0<r^m_{\phi(t)}\leqslant r^m_{\phi(0)}.
\end{equation}
In particular, the function $t\mapsto r^m_{\phi(t)}$ is decreasing.

Applying the maximum principle to (\ref{Revolution}), we have the following proposition.

\begin{proposition}\label{prop3.4}
Along the flow (\ref{flow}), there holds
$$\inf_MR_{\phi(t)}^m\geqslant \min\left\{\inf_MR_{\phi(0)}^m,0\right\}.$$
\end{proposition}

The following corollary follows from
(\ref{flow}) and Proposition \ref{prop3.4} immediately.

\begin{corollary}\label{cor3.5}
Along the flow (\ref{flow}), there holds
$$\frac{\partial}{\partial t}\phi(t)\geqslant \frac{m}{2}\left(\inf_MR_{\phi(0)}^m-r_{\phi(0)}^m\right).$$
\end{corollary}

For abbreviation, we let
\begin{equation}\label{sigma}
\sigma=\max\left\{\sup_M\big(1-R_{\phi(0)}^m\big),1\right\}.
\end{equation}
By Proposition \ref{prop3.4}, we have $R_{\phi(t)}^m+\sigma\geqslant 1$ for all $t\geqslant 0$.

\begin{lemma}\label{lem3.7}
For every $p>2$, we have
\begin{equation*}
\begin{split}
&\frac{d}{dt}\int_M(R_{\phi(t)}^m+\sigma)^{p-1}e^{-\phi(t)}dV_{g(t)}\\
&=-\frac{4(n+m-1)(p-2)}{p-1}\int_M|\nabla_{g(t)}(R_{\phi(t)}^m+\sigma)^{\frac{p-1}{2}}|^2
e^{-\phi(t)}dV_{g(t)}\\
&\hspace{4mm}
-\left(\frac{n+m+2}{2}-p\right)
\int_M\big((R_{\phi(t)}^m+\sigma)^{p-1}-(r_{\phi(t)}^m+\sigma)^{p-1}\big)
(R_{\phi(t)}^m-r_{\phi(t)}^m)e^{-\phi(t)}dV_{g(t)}\\
&\hspace{4mm}-(p-1)\sigma\int_M
\big((R_{\phi(t)}^m+\sigma)^{p-2}-(r_{\phi(t)}^m+\sigma)^{p-2}\big)
(R_{\phi(t)}^m-r_{\phi(t)}^m)e^{-\phi(t)}dV_{g(t)}.
\end{split}
\end{equation*}
\end{lemma}
\begin{proof}
This follows from differentiating
$\displaystyle\int_M(R_{\phi(t)}^m+\sigma)^{p-1}e^{-\phi(t)}dV_{g(t)}$
with respect to $t$, using the evolution equations
(\ref{Revolution}) and the integration by parts in (\ref{parts}).
This was done in \cite[Lemma 2.6]{Yan}.
The only difference between our case and \cite[Lemma 2.6]{Yan}
is that we have to take care of the term on the boundary $\partial M$.
But the term on the boundary $\partial M$ vanishes in view of (\ref{3.8}).
We leave the details to the reader.
\end{proof}

\begin{lemma}\label{lem3.8}
For every $p>\max\displaystyle\Big\{\frac{n+m}{2},2\Big\}$, we have
\begin{equation*}
\begin{split}
&\frac{d}{dt}\int_M|R_{\phi(t)}^m-r_{\phi(t)}^m|^pe^{-\phi(t)}dV_{g(t)}\leqslant C\int_M|R_{\phi(t)}^m-r_{\phi(t)}^m|^pe^{-\phi(t)}dV_{g(t)}\\
& +C\left(\int_M|R_{\phi(t)}^m-r_{\phi(t)}^m|^pe^{-\phi(t)}dV_{g(t)}\right)^{\frac{2p-(n+m)+2}{2p-(n+m)}}
\end{split}
\end{equation*}
for some uniform constant $C$ independent of $t$.
\end{lemma}
\begin{proof}
By (\ref{flow}) and (\ref{Revolution})
we compute
\begin{equation*}
\begin{split}
\frac{d}{dt}&\int_M|R_{\phi(t)}^m-r_{\phi(t)}^m|^pe^{-\phi(t)}dV_{g(t)}=p(n+m-1)\cdot\\
&\int_M|R_{\phi(t)}^m-r_{\phi(t)}^m|^{p-2}(R_{\phi(t)}^m-r_{\phi(t)}^m)
\big(\Delta_{\phi(t)}R_{\phi(t)}^m+R_{\phi(t)}^m(R_{\phi(t)}^m-r_{\phi(t)}^m)\big)e^{-\phi(t)}dV_{g(t)}\\
&\hspace{4mm}-\frac{n+m}{2}\int_M|R_{\phi(t)}^m-r_{\phi(t)}^m|^p(R_{\phi(t)}^m-r_{\phi(t)}^m)e^{-\phi(t)}dV_{g(t)}\\
&\hspace{4mm}
+\frac{(n+m-2)p}{2}\int_M|R_{\phi(t)}^m-r_{\phi(t)}^m|^{p-2}(R_{\phi(t)}^m-r_{\phi(t)}^m)
e^{-\phi(t)}dV_{g(t)}\\
&\hspace{8mm}\times\int_M(R_{\phi(t)}^m-r_{\phi(t)}^m)^2
e^{-\phi(t)}dV_{g(t)}.
\end{split}
\end{equation*}
Moreover, we have
\begin{equation*}
    \begin{split}
&\frac{d}{dt}\int_M|R_{\phi(t)}^m-r_{\phi(t)}^m|^pe^{-\phi(t)}dV_{g(t)}\\
&=-\frac{4(p-1)(n+m-1)}{p}
\int_M(R_{\phi(t)}^m-r_{\phi(t)}^m)^{\frac{p}{2}}L_{\phi(t)}^m\big((R_{\phi(t)}^m-r_{\phi(t)}^m)^{\frac{p}{2}}\big)e^{-\phi(t)}dV_{g(t)}\\
&\hspace{4mm}+
\left(\frac{(n+m-2)(p-1)}{p}+p-\frac{n+m}{2}\right)\int_M|R_{\phi(t)}^m-r_{\phi(t)}^m|^{p}(R_{\phi(t)}^m-r_{\phi(t)}^m)
e^{-\phi(t)}dV_{g(t)}\\
&\hspace{4mm}+\left(\frac{(n+m-2)(p-1)}{p}+p\right)r_{\phi(t)}^m\int_M|R_{\phi(t)}^m-r_{\phi(t)}^m|^{p}e^{-\phi(t)}dV_{g(t)}\\
&\hspace{4mm}
+\frac{(n+m-2)p}{2}\int_M|R_{\phi(t)}^m-r_{\phi(t)}^m|^{p-2}(R_{\phi(t)}^m-r_{\phi(t)}^m)
e^{-\phi(t)}dV_{g(t)}\\
&\hspace{8mm}\times\int_M(R_{\phi(t)}^m-r_{\phi(t)}^m)^2
e^{-\phi(t)}dV_{g(t)},
\end{split}
\end{equation*}
where we use (\ref{3.8}) and (\ref{parts})
in the last equality as follows:
\begin{equation*}
\begin{split}
&\int_M|R_{\phi(t)}^m-r_{\phi(t)}^m|^{p-2}(R_{\phi(t)}^m-r_{\phi(t)}^m)\Delta_{\phi(t)}R_{\phi(t)}^me^{-\phi(t)}dV_{g(t)}\\
&=-\int_M\langle\nabla_{g(t)}\big(|R_{\phi(t)}^m-r_{\phi(t)}^m|^{p-2}(R_{\phi(t)}^m-r_{\phi(t)}^m)\big),\nabla_{g(t)}R_{\phi(t)}^m\rangle e^{-\phi(t)}dV_{g(t)}\\
&\hspace{4mm}+\int_{\partial M}|R_{\phi(t)}^m-r_{\phi(t)}^m|^{p-2}(R_{\phi(t)}^m-r_{\phi(t)}^m)\frac{\partial R_{\phi(t)}^m}{\partial\nu_{g(t)}}e^{-\phi(t)}dA_{g(t)}\\
&=-(p-1)\int_M(R_{\phi(t)}^m-r_{\phi(t)}^m)^{p-2}|\nabla_{g(t)}R_{\phi(t)}^m|^2 e^{-\phi(t)}dV_{g(t)}\\
&=-\frac{4(p-1)}{p^2}
\int_M|\nabla_{g(t)}(R_{\phi(t)}^m-r_{\phi(t)}^m)^{\frac{p}{2}}|^2 e^{-\phi(t)}dV_{g(t)}\\
&\hspace{4mm}+\frac{4(p-1)}{p^2}
\int_{\partial M}(R_{\phi(t)}^m-r_{\phi(t)}^m)^{\frac{p}{2}}
\frac{\partial}{\partial\nu_{g(t)}}\big(
(R_{\phi(t)}^m-r_{\phi(t)}^m)^{\frac{p}{2}}\big) e^{-\phi(t)}dA_{g(t)}\\
&=\frac{4(p-1)}{p^2}\int_M
(R_{\phi(t)}^m-r_{\phi(t)}^m)^{\frac{p}{2}}
\Delta_{\phi(t)}\big(
(R_{\phi(t)}^m-r_{\phi(t)}^m)^{\frac{p}{2}}\big)e^{-\phi(t)}dV_{g(t)}.
\end{split}
\end{equation*}
Since $Y_{n,m}[(g_0, \phi_0)]> 0$
by assumption and the function $t\mapsto r_{\phi(t)}^m$
is monotonic decreasing, we have
\begin{equation*}
\begin{split}
&\frac{d}{dt}\int_M|R_{\phi(t)}^m-r_{\phi(t)}^m|^pe^{-\phi(t)}dV_{g(t)}\\
&\leq -\frac{(n+m-2)(p-1)}{p}
Y_{n,m}[(g_0, \phi_0)]
\left(\int_M|R_{\phi(t)}^m-r_{\phi(t)}^m|^{\frac{p(n+m)}{n+m-2}}
e^{-\phi(t)}dV_{g(t)}\right)^{\frac{n+m-2}{n+m}}\\
&\hspace{4mm}+
\left(\frac{(n+m-2)(p-1)}{p}+p-\frac{n+m}{2}\right)\int_M|R_{\phi(t)}^m-r_{\phi(t)}^m|^{p+1}e^{-\phi(t)}dV_{g(t)}\\
&\hspace{4mm}+\left(\frac{(n+m-2)(p-1)}{p}+p\right)r_{\phi(0)}^m\int_M|R_{\phi(t)}^m-r_{\phi(t)}^m|^{p}e^{-\phi(t)}dV_{g(t)}\\
&\hspace{4mm}
+\frac{(n+m-2)p}{2}\int_M|R_{\phi(t)}^m-r_{\phi(t)}^m|^{p-2}(R_{\phi(t)}^m-r_{\phi(t)}^m)
e^{-\phi(t)}dV_{g(t)}\\
&\hspace{8mm}\times\int_M(R_{\phi(t)}^m-r_{\phi(t)}^m)^2
e^{-\phi(t)}dV_{g(t)}.
\end{split}
\end{equation*}
By H\"{o}lder's inequality in $L^{p}(M, e^{-\phi(t)}dV_{g(t)})$ and (\ref{vol}), we have
\begin{equation*}
\begin{split}
&\int_M|R_{\phi(t)}^m-r_{\phi(t)}^m|^{p-2}(R_{\phi(t)}^m-r_{\phi(t)}^m)
e^{-\phi(t)}dV_{g(t)}\times\int_M(R_{\phi(t)}^m-r_{\phi(t)}^m)^2
e^{-\phi(t)}dV_{g(t)}\\
&\leq\left(\int_M|R_{\phi(t)}^m-r_{\phi(t)}^m|^{p}
e^{-\phi(t)}dV_{g(t)}\right)^{\frac{p+1}{p}}
\end{split}
\end{equation*}
and
\begin{equation*}
\begin{split}
\int_M|R_{\phi(t)}^m-r_{\phi(t)}^m|^{p+1}
e^{-\phi(t)}dV_{g(t)}&\leq\left(\int_M|R_{\phi(t)}^m-r_{\phi(t)}^m|^{p}
e^{-\phi(t)}dV_{g(t)}\right)^{\frac{2p-(n+m)+2}{2p}}\\
& \times\left(\int_M|R_{\phi(t)}^m-r_{\phi(t)}^m|^{\frac{p(n+m)}{n+m-2}}
e^{-\phi(t)}dV_{g(t)}\right)^{\frac{n+m-2}{2p}}.
\end{split}
\end{equation*}
Moreover, for any $\epsilon>0$,
we can apply the Young's inequality
to the last inequality
to deduce that
\begin{equation*}
\begin{split}
\int_M|R_{\phi(t)}^m-r_{\phi(t)}^m|^{p+1}
e^{-\phi(t)}dV_{g(t)}&\leq C_1\left(\int_M|R_{\phi(t)}^m-r_{\phi(t)}^m|^{p}
e^{-\phi(t)}dV_{g(t)}\right)^{\frac{2p-(n+m)+2}{2p-(n+m)}}\\
&+\epsilon
\left(\int_M|R_{\phi(t)}^m-r_{\phi(t)}^m|^{\frac{p(n+m)}{n+m-2}}
e^{-\phi(t)}dV_{g(t)}\right)^{\frac{n+m-2}{n+m}}
\end{split}
\end{equation*}
for some constant $C_1$.
Now the assertion follows from combining all these.
\end{proof}

In order to bound the solution $w(t)$ above and below in the interval $[0,T]$, we need the following lemmas.

\begin{lemma}\label{lem2.7}
Let $P$ be a smooth function on $(M, g, e^{-\phi}dV_g,e^{-\phi}dA_g, m)$. Moreover, assume that $w$ is a positive function such that
\begin{equation*}
-\frac{4(n+m-1)}{n+m-2}\Delta_{\phi}w+Pw\geqslant 0\mbox{ in }M
~~\mbox{ and }~~
\frac{\partial w}{\partial\nu_g}=0\mbox{ on }\partial M.
\end{equation*}
If $1\leqslant p<\frac{n}{n-2}$, there exists $C=C(n,m,p,g,\phi)$ and $r_0=r_0(M,g,\phi)$ such that
\begin{displaymath}
r^{-\frac{n}{p}}\|w\|_{L^p(B^+_{2r}(x))}\leqslant C\inf_{{B^+_{r}(x)}}w
\end{displaymath}
for any $x\in \partial M$, $r<r_0$ and $B^+_{r}(x)=M\cap B_{r}(x)$.
\end{lemma}
\begin{proof}
Without loss of generality, we may assume $r=1$. Let $\beta < 0$ and $0\leqslant \chi\in C^1_c(B^+_4)$. By assumption and integration by parts, we have
\begin{displaymath}
\int_M \langle dw, d(\chi^2 w^\beta)\rangle_g e^{-\phi}dV_g+C\int_M Pw^{\beta+1}\chi^2 e^{-\phi}dV_g \geqslant 0.
\end{displaymath}
Since $\beta<0$, we obtain
\begin{displaymath}
\begin{split}
|\beta|\int_M \chi^2 w^{\beta-1}|dw|^2_ge^{-\phi}dV_g &\leqslant 2C\int_M \chi w^\beta |d\chi|_g|dw|_g e^{-\phi}dV_g\\
&+ C\int_M |P|w^{\beta+1}\chi^2 e^{-\phi}dV_g.
\end{split}
\end{displaymath}
Applying Young's inequality to the first term on the right hand side, we arrive at
\begin{equation}\label{2.12}
\begin{split}
\int_M \chi^2 w^{\beta-1}|dw|^2_ge^{-\phi}dV_g&\leqslant C |\beta|^{-2} \int_M |d\chi|^2_g w^{\beta+1}e^{-\phi}dV_g\\
&+ C|\beta|^{-1}  \int_M |P|w^{\beta+1}\chi^2 e^{-\phi}dV_g.
\end{split}
\end{equation}
We set $u=w^{\frac{\beta+1}{2}}$, $\beta\neq -1$ such that (\ref{2.12}) can be rewritten as
\begin{equation}\label{2.13}
    \int_M \chi^2 |du|^2_ge^{-\phi}dV_g \leqslant C \int_M |d\chi|^2_g u^2 e^{-\phi}dV_g+ C \int_M |P|\chi^2 u^2 e^{-\phi}dV_g.
\end{equation}
In order to handle the right hand side of (\ref{2.13}), we use H\"older's and interpolation inequalities to get
\begin{equation}\label{2.14}
    \begin{split}
        \int_M |P|\chi^2 u^2 e^{-\phi}dV_g &\leqslant \|P\|_{L^{\frac{q}{2}}(B^+_4)}\|\chi u\|^2_{L^{\frac{2q}{q-2}}(B^+_4)}\\
        &\leqslant \|P\|_{L^{\frac{q}{2}}(B^+_4)} \left( \epsilon^{\frac{1}{2}}\|\chi u\|^2_{L^{\frac{2n}{n-2}}(B^+_4)}+\epsilon^{-\frac{\mu_1}{2}}\|\chi u\|^2_{L^{2}(B^+_4)} \right)\\
        &\leqslant \|P\|_{L^{\frac{q}{2}}(B^+_4)} \left( \epsilon\|\chi u\|^2_{L^{\frac{2n}{n-2}}(B^+_4)}+\epsilon^{-\mu_1}\|\chi u\|^2_{L^{2}(B^+_4)} \right),
    \end{split}
\end{equation}
where $\mu_1=\frac{n}{q-n}>0$. Choosing $\epsilon$ sufficiently small, we can make use of (\ref{2.12}), (\ref{2.13}) and (\ref{2.14}) to obtain
\begin{equation}\label{2.15}
    \left(\int_{B^+_4} (\chi u)^{\frac{2n}{n-2}}e^{-\phi}dV_g\right)^{\frac{n-2}{n}}\leqslant C(1+|\gamma|)^{2\mu_1+2}\int_{B^+_4} \left(|d\chi|^2_g+\chi^2\right)u^2 e^{-\phi}dV_g,
\end{equation}
where $\gamma=\beta+1<0$.

For any $1\leqslant r_a<r_b\leqslant 3$, we choose $\chi$ as a cut-off function satisfying $0\leqslant \chi \leqslant 1$, $|d\chi|\leqslant \frac{2}{r_b-r_a}$ and
\begin{displaymath}
\left\{
\begin{array}{ll}
     & \chi=1 ~~\mbox{ in }~~ B^+_{r_a},  \\
     & \chi=0 ~~\mbox{ in }~~ B^+_4\backslash B^+_{r_b}.
\end{array}
\right.
\end{displaymath}
Using this in (\ref{2.15}) yields
\begin{equation}
    \left(\int_{B^+_{r_a}} w^{\frac{\gamma n}{n-2}}e^{-\phi}dV_g\right)^{\frac{n-2}{n}}\leqslant \frac{C(1+|\gamma|)^{2\mu_1+2}}{r_b-r_a}\int_{B^+_{r_b}} w^{\gamma} e^{-\phi}dV_g.
\end{equation}
If we set $\Gamma(l,r)=\left(\int_{B^+_{r}} w^l e^{-\phi}dV_g\right)^{\frac{1}{l}}$ and $\delta=\frac{n}{n-2}$, the estimate above becomes
\begin{equation}\label{2.16}
    \Gamma(\gamma, r_b)\leqslant \left(\frac{C(1+|\gamma|)^{2\mu_1+2}}{r_b-r_a}\right)^{\frac{2}{|\gamma|}}\Gamma (\delta \gamma, r_a).
\end{equation}
It is well known that
\begin{displaymath}
\begin{split}
    \lim_{l\to +\infty} \Gamma(l,r)&=\sup_{B^+_r}w,\\
    \lim_{l\to -\infty} \Gamma(l,r)&=\inf_{B^+_r}w.
\end{split}
\end{displaymath}
The rest of the proof follows as in \cite{DN01} by iterating the estimate in (\ref{2.16}).
\end{proof}

\begin{lemma}\label{lem3.6}
Let $P$ be a smooth function on $(M, g, e^{-\phi}dV_g,e^{-\phi}dA_g, m)$.
Moreover, assume that $w$ is a positive function such that
\begin{equation*}
-\frac{4(n+m-1)}{n+m-2}\Delta_{\phi}w+Pw\geq 0\mbox{ in }M
~~\mbox{ and }~~
\frac{\partial w}{\partial\nu_g}=0\mbox{ on }\partial M.
\end{equation*}
There exists a constant $C$ depending only on  $(M, g, e^{-\phi}dV_g,e^{-\phi}dA_g, m)$
and $P$ such that
\begin{equation}\label{3.12}
\int_M we^{-\phi}dV_g\leq C\inf_M w.
\end{equation}
In particular, we have
\begin{equation}\label{3.13}
\int_M w^{\frac{2(n+m)}{n+m-2}}e^{-\phi}dV_g\leq C\Big(\inf_M w\Big)\Big(\sup_M w\Big)^{\frac{n+m+2}{n+m-2}}.
\end{equation}
\end{lemma}
\begin{proof}
Fix $r>0$ sufficiently small. Notice that the weighted Laplacian $\Delta_{\phi}$ has the same second-order terms as the classical Laplacian. The difference only occurs on lower order terms. Therefore, the interior weak Harnack inequality for linear elliptic equations \cite[Theorem 8.18]{DN01} can still hold in the weighted case, i.e. we obtain
\begin{equation}\label{2.20}
\int_{B_{2r}(x)} w e^{-\phi}dV_g\leqslant e^{-\inf\phi}\int_{B_{2r}(x)} w dV_g \leqslant e^{-\inf\phi}L_0\inf_{B_r(x)} w
\end{equation}
for some constant $L_0$, where $x\in M$ and $B_{2r}(x)\subset M$.

Combining (\ref{2.20}) with Lemma \ref{lem2.7} yields the global estimate
\begin{displaymath}
\int_{B^+_{2r}(x)} w e^{-\phi}dV_g \leqslant C\inf_{B^+_r(x)} w
\end{displaymath}
for some positive constant $C$, where $x\in M\cup \partial M$ and $B^+_{r}(x)=M\cap B_{r}(x)$. The assertion follows from the same argument as that in \cite[Proposition A.2]{Brendle4}.
\end{proof}

\begin{proposition}\label{prop3.8}
Given any $T>0$, we can find positive constants $C(T)$
and $c(T)$ such that
\begin{equation*}
c(T)\leq \inf_M w(t)\leq \sup_M w(t)\leq C(T)
\end{equation*}
for all $0\leq t\leq T$.
\end{proposition}
\begin{proof}
By Proposition \ref{prop3.4} and (\ref{boundr}), the conformal factor $w(t)$ satisfies
\begin{equation*}
\frac{\partial}{\partial t}w(t)
=-\frac{m+n-2}{4}(R^m_{\phi(t)}-r^m_{\phi(t)})w(t)
\leq \frac{m+n-2}{4}(r^m_{\phi(0)}+\sigma)w(t)
~~\mbox{ in }M.
\end{equation*}
Hence,
$$\frac{\partial}{\partial t}\log w(t)\leq\frac{m+n-2}{4}(r^m_{\phi(0)}+\sigma).$$
We conclude that $\displaystyle\sup_M w(t)\leq C(T)$ for all $0\leq t\leq T$.
Hence, if we define
$$P=R_{\phi_0}^m+\sigma\Big(\sup_{0\leq t\leq T}\sup_M w(t)\Big)^{\frac{4}{n+m-2}},$$
then we have
\begin{equation}\label{2.21}
\begin{split}
&-\frac{4(n+m-1)}{n+m-2}\Delta_{\phi_0}w(t)+Pw(t)\\
&\geq-\frac{4(n+m-1)}{n+m-2}\Delta_{\phi_0}w(t)+R_{\phi_0}^m w(t)+\sigma w(t)^{\frac{n+m+2}{n+m-2}}\\
&=(R_{\phi(t)}^m+\sigma)w(t)^{\frac{n+m+2}{n+m-2}}\geq 0
\end{split}
\end{equation}
for all $0\leq t\leq T$.
By (\ref{3.6}) and (\ref{2.21}), we can apply Lemma \ref{lem3.6} and find a positive constant $c(T)$
such that
\begin{displaymath}
\inf_M w(t)\Big(\sup_M w(t)\Big)^{\frac{n+m+2}{n+m-2}}
\geq c(T)
\end{displaymath}
for all $0\leq t\leq T$.
Since $\displaystyle\sup_M w(t)\leq C(T)$,
the assertion follows.
\end{proof}

\begin{proposition}
Let $0<\alpha<\displaystyle\frac{2m}{n+m}$.
Given any $T>0$, there exists a constant $C(T)$ such that
$$|w(x_1,t_1)-w(x_2,t_2)|\leq C(T)\big((t_1-t_2)^{\frac{\alpha}{2}}+
d(x_1,x_2)^\alpha\big)$$
for all $x_1, x_2\in M$
and $t_1,t_2\in [0,T]$ satisfying $0<t_1-t_2<1$.
Here, $d(x_1,x_2)$ is the distance between $x_1$ and $x_2$
with respect to the metric $g_0$.
\end{proposition}
\begin{proof}
By Lemma \ref{lem3.7} with $p=\displaystyle\frac{n+m+2}{2}$,
we obtain for all $0\leqslant t\leqslant T$
\begin{displaymath}
\frac{d}{dt}\int_M(R_{\phi(t)}^m+\sigma)^{\frac{n+m}{2}}e^{-\phi(t)}dV_{g(t)}
\leq 0,
\end{displaymath}
which implies for all $0\leqslant t\leqslant T$
\begin{displaymath}
\int_M(R_{\phi(t)}^m+\sigma)^{\frac{n+m}{2}}e^{-\phi(t)}dV_{g(t)}\leq C.
\end{displaymath}
This together with H\"{o}lder's inequality and (\ref{boundr})
implies that
\begin{equation}\label{3.14}
\begin{split}
&\left(\int_M|R_{\phi(t)}^m-r_{\phi(t)}^m|^{\frac{n+m}{2}}e^{-\phi(t)}dV_{g(t)}\right)^{\frac{2}{n+m}}\\
&\leq
\left(\int_M(R_{\phi(t)}^m+\sigma)^{\frac{n+m}{2}}e^{-\phi(t)}dV_{g(t)}\right)^{\frac{2}{n+m}}
+(r_{\phi(t)}^m+\sigma)\\
&\leq C.
\end{split}
\end{equation}
Let $\alpha=2-\displaystyle\frac{n}{p}$,
where $\displaystyle\frac{n}{2}<p<\frac{n+m}{2}$ with $m>0$.
Using (\ref{3.5}) and (\ref{3.14}) and Proposition \ref{prop3.8},
we obtain
\begin{equation}\label{3.15}
\int_M\left|\frac{4(n+m-1)}{n+m-2}\Delta_{\phi_0}w(t)
+R_{\phi_0}^mw(t)\right|^pe^{-\phi_0}dV_{g_0}\leq C(T)
\end{equation}
and
\begin{equation}\label{3.16}
\int_M\left|\frac{\partial}{\partial t}w(t)\right|^pe^{-\phi_0}dV_{g_0}\leq C(T)
\end{equation}
for all $t\in [0,T]$.
By the Sobolev embedding $W^{2,p}(M)\hookrightarrow C^{0,\alpha}(M)$,
the inequality (\ref{3.15}) implies that
$$|w(x_1,t)-w(x_2,t)|\leq C(T)
d(x_1,x_2)^\alpha$$
for all $x_1,x_2\in M$ and all $t\in [0,T]$.
Using (\ref{3.16}), we find
\begin{equation*}
\begin{split}
&|w(x,t_1)-w(x,t_2)|\\
&\leqslant C(t_1-t_2)^{-\frac{n}{2}}\int_{B_{\sqrt{t_1-t_2}(x)}}|w(x,t_1)-w(x,t_2)|e^{-\phi_0}dV_{g_0}\\
&\leqslant  C(t_1-t_2)^{-\frac{n}{2}}\int_{B_{\sqrt{t_1-t_2}(x)}}|w(t_1)-w(t_2)|e^{-\phi_0}dV_{g_0}
+C(T)(t_1-t_2)^{\frac{\alpha}{2}}\\
&\leqslant C(t_1-t_2)^{-\frac{n-2}{2}}\sup_{t_1\leqslant t\leqslant t_2}\int_{B_{\sqrt{t_1-t_2}(x)}}\left|\frac{\partial}{\partial t}w(t)\right|e^{-\phi_0}dV_{g_0}
+C(T)(t_1-t_2)^{\frac{\alpha}{2}}\\
&\leqslant C(t_1-t_2)^{\frac{\alpha}{2}}\sup_{t_1\leqslant t\leqslant t_2}\left(\int_{B_{\sqrt{t_1-t_2}(x)}}\left|\frac{\partial}{\partial t}w(t)\right|^pe^{-\phi_0}dV_{g_0}\right)^{\frac{1}{p}}
+C(T)(t_1-t_2)^{\frac{\alpha}{2}}\\
&\leqslant C(T)(t_1-t_2)^{\frac{\alpha}{2}},
\end{split}
\end{equation*}
for all $x\in M$ and all $t_1,t_2\in [0,T]$ satisfying $0<t_1-t_2<1$.
This proves the assertion.
\end{proof}

Now we can use the standard regularity theory for
parabolic equations
to show that all higher order derivatives of $w(t)$
are uniformly bounded on every fixed time interval $[0,T]$.
Therefore, the flow exists for all time.

\section{Proof of the main result assuming Proposition \ref{prop4.3}}

In this section, we will prove Theorem \ref{main}
by assuming
Proposition \ref{prop4.3}.
In the following, $c$ and $C$ are positive constants
independent of $t$, and may change from line to line.

\begin{proposition}\label{prop4.1}
For any $\displaystyle\max\Big\{\frac{n+m}{2},2\Big\}<p<
\frac{n+m+2}{2}$, we have
$$\lim_{t\to\infty}\int_M|R_{\phi(t)}^m-r_{\phi(t)}^m|^pe^{-\phi(t)}dV_{g(t)}
=0.$$
\end{proposition}
\begin{proof}
Since $p>2$, it follows from Lemma  \ref{lem3.7} that
\begin{equation*}
\begin{split}
&\frac{d}{dt}\int_M(R_{\phi(t)}^m+\sigma)^{p-1}e^{-\phi(t)}dV_{g(t)}\\
&\leqslant
-\left(\frac{n+m+2}{2}-p\right)
\int_M\big((R_{\phi(t)}^m+\sigma)^{p-1}-(r_{\phi(t)}^m+\sigma)^{p-1}\big)
(R_{\phi(t)}^m-r_{\phi(t)}^m)e^{-\phi(t)}dV_{g(t)}.
\end{split}
\end{equation*}
Since $p>2$, we have
$$\big((R_{\phi(t)}^m+\sigma)^{p-1}-(r_{\phi(t)}^m+\sigma)^{p-1}\big)
(R_{\phi(t)}^m-r_{\phi(t)}^m)
\geqslant c|R_{\phi(t)}^m-r_{\phi(t)}^m|^p$$
for some constant $c>0$.
Since $p<\displaystyle\frac{n+m+2}{2}$,
we obtain
$$
\frac{d}{dt}\int_M(R_{\phi(t)}^m+\sigma)^{p-1}e^{-\phi(t)}dV_{g(t)}
\leqslant -c\int_M|R_{\phi(t)}^m-r_{\phi(t)}^m|^pe^{-\phi(t)}dV_{g(t)}.$$
Integrating it with respect to $t$ yields
$$\int_0^\infty\int_M|R_{\phi(t)}^m-r_{\phi(t)}^m|^pe^{-\phi(t)}dV_{g(t)}dt\leqslant C.$$
In particular, we have
$$\liminf_{t\to\infty}\int_M|R_{\phi(t)}^m-r_{\phi(t)}^m|^pe^{-\phi(t)}dV_{g(t)}
=0.$$
On the other hand, since $p>\displaystyle\max\Big\{\frac{n+m}{2},2\Big\}$,
it follows from Lemma \ref{lem3.8}
that
\begin{equation*}
\begin{split}
&\frac{d}{dt}\int_M|R_{\phi(t)}^m-r_{\phi(t)}^m|^pe^{-\phi(t)}dV_{g(t)}\\
&\leqslant C\left(\int_M|R_{\phi(t)}^m-r_{\phi(t)}^m|^pe^{-\phi(t)}dV_{g(t)}\right)^{\frac{2p-(n+m)+2}{2p-(n+m)}}
+C\int_M|R_{\phi(t)}^m-r_{\phi(t)}^m|^pe^{-\phi(t)}dV_{g(t)}.
\end{split}
\end{equation*}
From this, the assertion follows.
\end{proof}

Hence, if we define
\begin{equation}\label{4.1}
r_\infty^m=\lim_{t\to\infty} r_{\phi(t)}^m,
\end{equation}
then, we have the following result:

\begin{corollary}\label{cor4.2}
For every $1<p<\displaystyle\frac{n+m+2}{2}$, we have
$$\lim_{t\to\infty} \int_M|R_{\phi(t)}^m-r_\infty^m|^pe^{-\phi(t)}dV_{g(t)}
=0.$$
\end{corollary}
\begin{proof}
It follows from H\"{o}lder's inequality and Proposition \ref{prop4.1} that
\begin{equation}\label{4.2}
\lim_{t\to\infty} \int_M|R_{\phi(t)}^m-r_{\phi(t)}^m|^pe^{-\phi(t)}dV_{g(t)}
=0
\end{equation}
for all  $1<p<\displaystyle\frac{n+m+2}{2}$.
By  Minkowski inequality, we have
\begin{equation*}
\begin{split}
&\left(\int_M|R_{\phi(t)}^m-r_\infty^m|^pe^{-\phi(t)}dV_{g(t)}\right)^{\frac{1}{p}}\\
&\leqslant
\left(\int_M|R_{\phi(t)}^m-r_{\phi(t)}^m|^pe^{-\phi(t)}dV_{g(t)}\right)^{\frac{1}{p}}
+\left(|r_{\phi(t)}^m-r_\infty^m|^p\int_Me^{-\phi(t)}dV_{g(t)}\right)^{\frac{1}{p}}.
\end{split}
\end{equation*}
Together with (\ref{4.1}) and (\ref{4.2}), this implies the assertion.
\end{proof}

The proof of Theorem \ref{main} will be based on the following proposition.

\begin{proposition}\label{prop4.3}
Let $\{t_i:i\in\mathbb{N}\}$ be a sequence of times such that
$t_i\to\infty$ as $i\to\infty$.
Then we can find a real number $0<\gamma<1$ and a constant $C$ such that,
after passing to a subsequence, we have
\begin{equation}\label{4.3}
r_{\phi(t_i)}^m-r_\infty^m
\leqslant C\left(\int_M
|R_{\phi(t_i)}^m-r_\infty^m|^{\frac{2(n+m)}{n+m+2}}w(t_i)^{\frac{2(n+m)}{n+m-2}}e^{-\phi_0}dV_{g_0}\right)^{\frac{n+m+2}{2(n+m)}(1+\gamma)}
\end{equation}
for all integers $i$ in that sequence.
Note that $\gamma$ and $C$ may depend
on the sequence $\{t_i: i\in\mathbb{N}\}$.
\end{proposition}

The following result is an immediate consequence of Proposition \ref{prop4.3}.

\begin{proposition}\label{prop4.4}
There exist real numbers $0<\gamma<1$
and $t_0>0$ such that
\begin{equation}\label{4.4}
r_{\phi(t)}^m-r_\infty^m
\leqslant C\left(\int_M
|R_{\phi(t)}^m-r_\infty^m|^{\frac{2(n+m)}{n+m+2}}w(t)^{\frac{2(n+m)}{n+m-2}}e^{-\phi_0}dV_{g_0}\right)^{\frac{n+m+2}{2(n+m)}(1+\gamma)}
\end{equation}
for all $t\geq t_0$.
\end{proposition}
\begin{proof}
Suppose this is not true.
Then there exists a sequence of times
$\{t_i: i\in\mathbb{N}\}$ such that
$t_i\geqslant i$ and
\begin{displaymath}
r_{\phi(t_i)}^m-r_\infty^m
\geqslant C\left(\int_M
|R_{\phi(t_i)}^m-r_\infty^m|^{\frac{2(n+m)}{n+m+2}}w(t_i)^{\frac{2(n+m)}{n+m-2}}e^{-\phi_0}dV_{g_0}\right)^{\frac{n+m+2}{2(n+m)}(1+\frac{1}{i})}.
\end{displaymath}
We now apply Proposition \ref{prop4.3}
to this sequence $\{t_i: i\in\mathbb{N}\}$.
Hence, there exist an infinite subset $I\subset\mathbb{N}$
and real numbers $0<\gamma<1$ and $C$ such that
\begin{displaymath}
r_{\phi(t_i)}^m-r_\infty^m \leqslant C\left(\int_M
|R_{\phi(t_i)}^m-r_\infty^m|^{\frac{2(n+m)}{n+m+2}}w(t_i)^{\frac{2(n+m)}{n+m-2}}e^{-\phi_0}dV_{g_0}\right)^{\frac{n+m+2}{2(n+m)}(1+\gamma)}
\end{displaymath}
for all $i\in I$. Thus, we conclude that
$$1\leqslant C\left(\int_M
|R_{\phi(t_i)}^m-r_\infty^m|^{\frac{2(n+m)}{n+m+2}}w(t_i)^{\frac{2(n+m)}{n+m-2}}e^{-\phi_0}dV_{g_0}\right)^{\frac{n+m+2}{2(n+m)}(\gamma-\frac{1}{i})}
$$
for all $i\in I$.

On the other hand, it follows from
Corollary \ref{cor3.5} with $p=\frac{2(n+m)}{n+m+2}<\frac{n+m+2}{2}$
and $ w(t_i)^{\frac{2(n+m)}{n+m-2}}e^{-\phi_0}dV_{g_0}=e^{-\phi(t_i)}dV_{g(t_i)}$
that
$$\lim_{i\to\infty}\int_M|R_{\phi(t_i)}^m-r_\infty^m|^{\frac{2(n+m)}{n+m+2}}
e^{-\phi(t_i)}dV_{g(t_i)}=0.$$
Therefore, if $i$ is sufficiently large,
\begin{equation*}
\begin{split}
1&\leqslant \lim_{i\to\infty}\left(\int_M
|R_{\phi(t_i)}^m-r_\infty^m|^{\frac{2(n+m)}{n+m+2}}w(t_i)^{\frac{2(n+m)}{n+m-2}}e^{-\phi_0}dV_{g_0}\right)^{\frac{n+m+2}{2(n+m)}(\gamma-\frac{1}{i})}\\
&\leqslant
\lim_{i\to\infty}\left(\int_M
|R_{\phi(t_i)}^m-r_\infty^m|^{\frac{2(n+m)}{n+m+2}}w(t_i)^{\frac{2(n+m)}{n+m-2}}e^{-\phi_0}dV_{g_0}\right)^{\frac{n+m+2}{2(n+m)}\frac{\gamma}{2}}
=0,
\end{split}
\end{equation*}
which is a contradiction.
\end{proof}

\begin{proposition}\label{prop4.5}
There holds
\begin{displaymath}
\int_0^\infty\left(
\int_M
|R_{\phi(t)}^m-r_{\phi(t)}^m|^2w(t)^{\frac{2(n+m)}{n+m-2}}e^{-\phi_0}dV_{g_0}\right)^{\frac{1}{2}}dt\leq C.
\end{displaymath}
\end{proposition}

\begin{proof}
It follows from Proposition  \ref{prop4.4}
that
\begin{equation*}
\begin{split}
r_{\phi(t)}^m-r_\infty^m
&\leqslant C\left(\int_M
|R_{\phi(t)}^m-r_{\phi(t)}^m|^{\frac{2(n+m)}{n+m+2}}w(t)^{\frac{2(n+m)}{n+m-2}}e^{-\phi_0}dV_{g_0}\right)^{\frac{n+m+2}{2(n+m)}(1+\gamma)}\\
&\hspace{4mm}+C\big(r_{\phi(t)}^m-r_\infty^m\big)^{1+\gamma}.
\end{split}
\end{equation*}
Hence, we have
$$r_{\phi(t)}^m-r_\infty^m
\leqslant C\left(\int_M
|R_{\phi(t)}^m-r_{\phi(t)}^m|^{\frac{2(n+m)}{n+m+2}}w(t)^{\frac{2(n+m)}{n+m-2}}e^{-\phi_0}dV_{g_0}\right)^{\frac{n+m+2}{2(n+m)}(1+\gamma)} $$
if $t$ is sufficiently large.
Therefore, by H\"{o}lder's inequality
and (\ref{vol}), we have
\begin{equation}\label{4.5}
\begin{split}
\frac{d}{dt}\big(r_{\phi(t)}^m-r_\infty^m\big)
&=-\frac{n+m-2}{2}\int_M(R_{\phi(t)}^m-r_{\phi(t)}^m)^2w(t)^{\frac{2(n+m)}{n+m-2}}e^{-\phi_0}dV_{g_0}\\
&\leqslant
-\frac{n+m-2}{2}\left(\int_M|R_{\phi(t)}^m-r_{\phi(t)}^m|^{\frac{2(n+m)}{n+m+2}}w(t)^{\frac{2(n+m)}{n+m-2}}e^{-\phi_0}dV_{g_0}\right)^{\frac{n+m+2}{n+m}}\\
&\leqslant
-\frac{n+m-2}{2}\big(r_{\phi(t)}^m-r_\infty^m\big)^{\frac{2}{1+\gamma}}.
\end{split}
\end{equation}
This implies that
$$\frac{d}{dt}\big(r_{\phi(t)}^m-r_\infty^m\big)^{-\frac{1-\gamma}{1+\gamma}}\geqslant c.$$
From this, we can deduce that if $t$ is sufficiently large
$$r_{\phi(t)}^m-r_\infty^m\leqslant Ct^{-\frac{1-\gamma}{1+\gamma}}.$$
Integrating the first equality in (\ref{4.5})
from $T$ to $2T$ yields
$$r_{\phi(T)}^m-r_{\phi(2T)}^m
=
\frac{n+m-2}{2}\int_T^{2T}\int_M(R_{\phi(t)}^m-r_{\phi(t)}^m)^2w(t)^{\frac{2(n+m)}{n+m-2}}
e^{-\phi_0}dV_{g_0}dt.$$
Hence, by H\"{o}lder's inequality, we find
\begin{equation*}
\begin{split}
&\int_T^{2T}\left(\int_M(R_{\phi(t)}^m-r_{\phi(t)}^m)^2w(t)^{\frac{2(n+m)}{n+m-2}}
e^{-\phi_0}dV_{g_0}\right)^{\frac{1}{2}}dt\\
&\leqslant
\left(T\int_T^{2T}\int_M(R_{\phi(t)}^m-r_{\phi(t)}^m)^2w(t)^{\frac{2(n+m)}{n+m-2}}
e^{-\phi_0}dV_{g_0}dt\right)^{\frac{1}{2}}\\
&\leqslant \left(\frac{2}{n+m-2}T\big(r_{\phi(T)}^m-r_{\phi(2T)}^m\big)\right)^{\frac{1}{2}}\\
&\leqslant CT^{-\frac{\gamma}{1+\gamma}}
\end{split}
\end{equation*}
if $T$ is sufficiently large. Since $0<\gamma<1$, we can conclude that
\begin{equation*}
\begin{split}
&\int_0^\infty\left(\int_M(R_{\phi(t)}^m-r_{\phi(t)}^m)^2w(t)^{\frac{2(n+m)}{n+m-2}}
e^{-\phi_0}dV_{g_0}\right)^{\frac{1}{2}}dt\\
&=\int_0^1\left(\int_M(R_{\phi(t)}^m-r_{\phi(t)}^m)^2w(t)^{\frac{2(n+m)}{n+m-2}}
e^{-\phi_0}dV_{g_0}\right)^{\frac{1}{2}}dt\\
&\hspace{4mm}+\sum_{k=0}^\infty\int_{2^k}^{2^{k+1}}\left(\int_M(R_{\phi(t)}^m-r_{\phi(t)}^m)^2w(t)^{\frac{2(n+m)}{n+m-2}}
e^{-\phi_0}dV_{g_0}\right)^{\frac{1}{2}}dt\\
&\leqslant C\left(1+\sum_{k=0}^\infty 2^{-\frac{\gamma}{1+\gamma}k}\right)\leqslant C.
\end{split}
\end{equation*}
This proves the assertion.
\end{proof}

\begin{proposition}\label{prop4.6}
Given any $\eta_0>0$, we can find a real number $r>0$ such that
\begin{equation}\label{4.6}
\int_{B^+_r(x)}w(t)^{\frac{2(n+m)}{n+m-2}}e^{-\phi_0}dV_{g_0}\leqslant \eta_0
\end{equation}
for all $x\in M\cup \partial M$ and $t\geqslant 0$, where $B^+_{r}(x)=M\cap B_{r}(x)$.
\end{proposition}
\begin{proof}
It follows from Proposition \ref{prop4.5} that
we can find a real number $T>0$ such that
\begin{equation}\label{4.7}
\int_0^\infty\left(
\int_M
|R_{\phi(t)}^m-r_{\phi(t)}^m|^2w(t)^{\frac{2(n+m)}{n+m-2}}e^{-\phi_0}dV_{g_0}\right)^{\frac{1}{2}}dt\leqslant
\frac{\eta_0}{4(n+m)}.
\end{equation}
By Proposition \ref{prop3.8}, we can choose a real number $r>0$ such that
\begin{equation}\label{4.8}
\int_{B^+_r(x)}w(t)^{\frac{2(n+m)}{n+m-2}}e^{-\phi_0}dV_{g_0}\leqslant \frac{\eta_0}{2}
\end{equation}
for all $x\in M\cup \partial M$ and $0\leqslant t\leqslant T$.
By (\ref{vol}) and H\"{o}lder's inequality, we have
\begin{equation*}
\begin{split}
\frac{d}{dt}\int_{B^+_r(x)}e^{-\phi(t)}dV_{g(t)}&=
-\frac{n+m}{2}\int_{B^+_r(x)}(R_{\phi(t)}^m-r_{\phi(t)}^m)e^{-\phi(t)}dV_{g(t)}\\
&\leqslant \frac{n+m}{2}\left(\int_{M}(R_{\phi(t)}^m-r_{\phi(t)}^m)^2e^{-\phi(t)}dV_{g(t)}\right)^{\frac{1}{2}}.
\end{split}
\end{equation*}
Integrating this over $[T,t]$ yields
\begin{equation*}
\begin{split}
&\int_{B^+_r(x)}w(t)^{\frac{2(n+m)}{n+m-2}}e^{-\phi_0}dV_{g_0}\\
&\leqslant \int_{B^+_r(x)}w(T)^{\frac{2(n+m)}{n+m-2}}e^{-\phi_0}dV_{g_0}
+\frac{n+m}{2}\int_T^\infty\left(\int_{M}(R_{\phi(t)}^m-r_{\phi(t)}^m)^2e^{-\phi(t)}dV_{g(t)}\right)^{\frac{1}{2}}dt\\
&\leqslant \eta_0.
\end{split}
\end{equation*}
for all $x\in M$ and all $t\geqslant T$,
where we have used (\ref{4.7}) and (\ref{4.8}) in the last inequality.
This proves the assertion.
\end{proof}

\begin{lemma}\label{lem4.7}
Let $p=\frac{2(n+m)}{n+m-2}$ and $q>\frac{n}{2}$.
There are positive constants $\eta_1$ and $C$ such that if
\begin{equation}
\begin{array}{cc}
&g=w^{\frac{4}{n+m-2}}g_0, \\
&e^{-\phi}=w^{\frac{2m}{n+m-2}}e^{-\phi_0},
\end{array}
\end{equation}
and
$$\int_{B^+_{4r}(x)}e^{-\phi}dV_g\leq 1~~\mbox{ and }~~\int_{B^+_{4r}(x)}|R^m_\phi|^q e^{-\phi}dV_g\leq \eta_1,$$
where $B^+_{4r}(x)=M\cap B_{4r}(x)$ is the geodesic ball with respect to $g_0$ and $r<1$, then
$$w(x)\leq Cr^{-\frac{n}{p}}\left(\int_{B^+_{4r}(x)}e^{-\phi}dV_g\right)^{\frac{1}{p}}.$$
\end{lemma}
\begin{proof}
By the smoothness of the conformal factor $w(t)$, there exists $r_0$ a real number such that $r_0<r$ and
\begin{displaymath}
(r-s)^{\frac{n}{p}}\sup_{B^+_s(x)} w\leqslant (r-r_0)^{\frac{n}{p}}\sup_{B^+_{r_0}(x)} w
\end{displaymath}
for all $s<r$. Moreover, we choose a point $x_0\in \overline{B^+_{r_0}(x)}$ such that
\begin{displaymath}
\sup_{B^+_{r_0}(x)} w=w(x_0).
\end{displaymath}
Notice that the conformal weighted Laplacian $L^m_{\phi_0}$ has the same leading term as the classical Laplacian $\Delta_{g_0}$. The difference only occurs on lower order terms. If $x_0$ is in the interior of $M$, using a standard interior estimate for linear elliptic equations in \cite[Theorem 8.17]{DN01}, we obtain
\begin{equation}
\begin{split}
s^{\frac{n}{p}}w(x_0)&\leqslant C\left(\int_{B_s(x_0)}w^{p}e^{-\phi_0}dV_{g_0}\right)^{\frac{1}{p}}\\
& +Cs^{\frac{n}{p}+2-\frac{n}{q}}\left(\int_{B_s(x_0)}\left|\frac{4(n+m-1)}{(n+m-2)}L^m_{\phi_0}w\right|^q e^{-\phi_0}dV_{g_0}\right)^{\frac{1}{q}}
\end{split}
\end{equation}
for $s\leqslant \frac{r-r_0}{2}$ and $B_s(x_0)\subset M$.

If $x_0$ is on the boundary $\partial M$, we may adapt the argument in \cite[Porposition A-2]{Almaraz&Sun} to obtain
\begin{equation}
\begin{split}
    s^{\frac{n}{p}}w(x_0)&\leqslant C\left(\int_{B^+_{2s}(x_0)}w^{p}e^{-\phi_0}dV_{g_0}\right)^{\frac{1}{p}}\\
& +Cs^{\frac{n}{p}+2-\frac{n}{q}}\left(\int_{B^+_{4s}(x_0)}\left|\frac{4(n+m-1)}{(n+m-2)}L^m_{\phi_0}w\right|^q e^{-\phi_0}dV_{g_0}\right)^{\frac{1}{q}},
\end{split}
\end{equation}
for $s<\tilde{r}$, where $\tilde{r}$ is the constant in \cite[Porposition A-2]{Almaraz&Sun}.

In both cases we have
\begin{displaymath}
\begin{split}
    s^{\frac{n}{p}}w(x_0)&\leqslant C\left(\int_{B^+_{4s}(x_0)}w^{p}e^{-\phi_0}dV_{g_0}\right)^{\frac{1}{p}}\\
& +Cs^{\frac{n}{p}+2-\frac{n}{q}}\left(\int_{B^+_{4s}(x_0)}\left|\frac{4(n+m-1)}{(n+m-2)}L^m_{\phi_0}w\right|^q e^{-\phi_0}dV_{g_0}\right)^{\frac{1}{q}},
\end{split}
\end{displaymath}
for $s<\min \{\frac{r-r_0}{2},\tilde{r}\}$. The assertion follows from the same iteration argument as that in \cite[Proposition A.2]{Brendle4}.
\end{proof}

\begin{proposition}\label{prop4.8}
Along the flow, the function $w(t)$ satisfies
\begin{equation*}
c\leqslant \inf_M w(t)\leqslant \sup_M w(t)\leqslant C
\end{equation*}
for all $t\geq 0$. Here, $c$ and $C$ are positive constants independent of $t$.
\end{proposition}
\begin{proof}
Fix $\displaystyle\frac{n}{2}<q<p<\frac{n+m+2}{2}$.
It follows from Corollary \ref{cor4.2}
that
$$\int_M|R_{\phi(t)}^m|^pe^{-\phi(t)}dV_{g(t)}\leq C,$$
for some constant $C$ independent of $t$.
By Proposition \ref{prop4.6}, we can find a constant
$r>0$ independent of $t$ such that
$$\int_{B^+_{4r}(x)}w(t)^{\frac{2(n+m)}{n+m-2}}e^{-\phi_0}dV_{g_0}\leqslant \eta_0$$
for all $x\in M\cup \partial M$ and $t\geq 0$.
By H\"{o}lder's inequality, we have
$$\int_{B^+_{4r}(x)}|R_{\phi(t)}^m|^qe^{-\phi(t)}dV_{g(t)}
\leqslant \left(\int_{B^+_{4r}(x)}e^{-\phi(t)}dV_{g(t)}
\right)^{\frac{p-q}{p}}
\left(\int_M|R_{\phi(t)}^m|^pe^{-\phi(t)}dV_{g(t)}\right)^{\frac{q}{p}}.
$$
Hence, if we choose $\eta_0$ sufficiently small, we then have
$$\int_{B^+_{4r}(x)}|R_{\phi(t)}^m|^qe^{-\phi(t)}dV_{g(t)}\leqslant \eta_1$$
for all $x\in M\cup \partial M$ and all $t\geq 0$. Here, $\eta_1$
is the constant appearing in Lemma \ref{lem4.7}. We can now apply Lemma \ref{lem4.7}
at the maximum point of $w(t)$ to deduce that
\begin{equation*}
\sup_M w(t)\leqslant Cr^{-\frac{n}{p}}\left(\int_{B^+_{4r}(x)}e^{-\phi(t)}dV_{g(t)}\right)^{\frac{1}{p}}.
\end{equation*}
Together with (\ref{vol}), this implies that $w(t)$ is uniformly bounded from above.
Hence, if we define
$$P=R_{\phi_0}^m+\sigma\Big(\sup_{t\geq 0}\sup_M w(t)\Big)^{\frac{4}{n+m-2}}$$
where $\sigma$ is given as in (\ref{sigma}),
then we have
\begin{equation*}
\begin{split}
&-\frac{4(n+m-1)}{n+m-2}\Delta_{\phi_0}w(t)+Pw(t)\\
&\geqslant -\frac{4(n+m-1)}{n+m-2}\Delta_{\phi_0}w(t)+R_{\phi_0}^mw(t)+\sigma w(t)^{\frac{n+m+2}{n+m-2}}\\
&=(R_{\phi(t)}^m+\sigma)w(t)^{\frac{n+m+2}{n+m-2}}\geq 0.
\end{split}
\end{equation*}
By (\ref{vol}) and Lemma \ref{lem3.6}, we can find a positive constant $c$ independent of $t$
such that
$$\inf_M w(t)\Big(\sup_M w(t)\Big)^{\frac{n+m+2}{n+m-2}}
\geqslant c
$$
for all $t\geq 0$.
This implies that $w(t)$ is uniformly bounded from below,
since $w(t)$ is uniformly bounded from above.
This proves the assertion.
\end{proof}

\begin{proposition}\label{Holder}
Let $0<\alpha<\displaystyle\frac{2m}{n+m}$. There holds
$$|w(x_1,t_1)-w(x_2,t_2)|\leqslant C\big((t_1-t_2)^{\frac{\alpha}{2}}+
d(x_1,x_2)^\alpha\big)$$
for all $x_1, x_2\in M$
and $0<t_1-t_2<1$.
Here, $C$ is a positive constant independent of $t_1$ and $t_2$.
\end{proposition}
\begin{proof}
Let $\alpha=2-\displaystyle\frac{n}{p}$, where
$\displaystyle\frac{n}{2}<p<\frac{n+m}{2}$.
As in the proof of Proposition \ref{prop3.8},
we can deduce from Proposition \ref{prop4.8} that
\begin{equation}\label{4.9}
\int_M\left|\frac{4(n+m-1)}{n+m-2}\Delta_{\phi_0}w(t)
+R_{\phi_0}^mw(t)\right|^pe^{-\phi_0}dV_{g_0}\leqslant C
\end{equation}
and
\begin{equation}\label{4.10}
\int_M\left|\frac{\partial}{\partial t}w(t)\right|^pe^{-\phi_0}dV_{g_0}\leq C
\end{equation}
where $C$ is a positive constant independent of $t$.
By the Sobolev embedding $W^{2,p}(M)\hookrightarrow C^{0,\alpha}(M)$,
the inequality (\ref{4.9}) implies that
$$|w(x_1,t)-w(x_2,t)|\leq C
d(x_1,x_2)^\alpha$$
for all $x_1,x_2\in M$ and all $t\geq 0$.
On the other hand, by the second inequality (\ref{4.10}),
we find \begin{equation*}
\begin{split}
&|w(x,t_1)-w(x,t_2)|\\
&\leq C(t_1-t_2)^{-\frac{n}{2}}\int_{B_{\sqrt{t_1-t_2}(x)}}|w(x,t_1)-w(x,t_2)|e^{-\phi_0}dV_{g_0}\\
&\leqslant  C(t_1-t_2)^{-\frac{n}{2}}\int_{B_{\sqrt{t_1-t_2}(x)}}|w(t_1)-w(t_2)|e^{-\phi_0}dV_{g_0}
+C(t_1-t_2)^{\frac{\alpha}{2}}\\
&\leqslant C(t_1-t_2)^{-\frac{n-2}{2}}\sup_{t_1\leqslant t\leqslant t_2}\int_{B_{\sqrt{t_1-t_2}(x)}}\left|\frac{\partial}{\partial t}w(t)\right|e^{-\phi_0}dV_{g_0}
+C(t_1-t_2)^{\frac{\alpha}{2}}\\
&\leqslant C(t_1-t_2)^{\frac{\alpha}{2}}\sup_{t_1\leqslant t\leqslant t_2}\left(\int_{B_{\sqrt{t_1-t_2}(x)}}\left|\frac{\partial}{\partial t}w(t)\right|^pe^{-\phi_0}dV_{g_0}\right)^{\frac{1}{p}}
+C(t_1-t_2)^{\frac{\alpha}{2}}\\
&\leqslant C(t_1-t_2)^{\frac{\alpha}{2}},
\end{split}
\end{equation*}
for all $x\in M$ and all $t_1,t_2\geqslant 0$ satisfying $0<t_1-t_2<1$.
This proves the assertion.
\end{proof}

Now we can use the standard regularity theory for
parabolic equations
to show that all higher order derivatives of $w(t)$
are uniformly bounded on $[0,\infty)$. The uniqueness of the asymptotic limit follows
Proposition \ref{prop4.5}.
This completes the proof of Theorem \ref{main}.

\section{Proof of Proposition \ref{prop4.3}}

Let $\{t_i, i\in\mathbb{N}\}$ be a sequence of times such that $t_i\to\infty$ as $i\to\infty$.
For abbreviation, we let $w_i=w(t_i)$.
The normalization condition (\ref{vol}) implies that
\begin{equation}\label{5.1}
\int_Mw_i^{\frac{2(n+m)}{n+m-2}}e^{-\phi_0}dV_{g_0}=1
\end{equation}
for all $i\in\mathbb{N}$.
Moreover, it follows from Corollary \ref{cor4.2} that
$$\int_M|R_{\phi(t_i)}^m-r_\infty^m|^{\frac{2(n+m)}{n+m+2}}e^{-\phi(t_i)}dV_{g(t_i)}\to0,$$
and hence
\begin{equation}\label{5.2}
\int_M\left|\frac{4(n+m-1)}{n+m-2}\Delta_{\phi_0}w_i-R_{\phi_0}^m w_i+r_\infty^mw_i^{\frac{n+m+2}{n+m-2}}\right|^{\frac{2(n+m)}{n+m+2}}e^{-\phi_0}dV_{g_0}\to0
\end{equation}
as $i\to\infty$. On the other hand, it follows from (\ref{3.6}) that
\begin{equation}\label{5.3}
\frac{\partial w_i}{\partial\nu_{g_0}}=0~~\mbox{ on }\partial M
\end{equation}
for all $i\in\mathbb{N}$.
By the standard elliptic theory, we have the following compactness result.

\begin{proposition}\label{prop5.1}
Let $\{w_i: i\in\mathbb{N}\}$ be a sequence of positive functions satisfying
\eqref{5.1} and \eqref{5.2}. After passing to a subsequence if necessary,
$w_i$ converges to a positive smooth function $w_\infty$
satisfying
\begin{equation*}
\frac{4(n+m-1)}{n+m-2}\Delta_{\phi_0}w_\infty-R_{\phi_0}^m w_\infty+r_\infty^mw_\infty^{\frac{n+m+2}{n+m-2}}=0\mbox{ in }M
~~\mbox{ and }~~
\frac{\partial w_\infty}{\partial\nu_{g_0}}=0\mbox{ on }\partial M.
\end{equation*}
\end{proposition}
\begin{proof}
Since $\frac{n+m+2}{n+m-2}<\frac{n+2}{n-2}$,
the assertion now follows from
(\ref{5.1})-(\ref{5.3})
and the standard elliptic theory
 \cite[Section 8, Theorem 3]{Evans10}.
\end{proof}

In order to prove Proposition \ref{prop4.3}, we need the following:

\begin{proposition}\label{prop5.2}
There exists a sequence of smooth functions $\{\psi_a: a\in\mathbb{N}\}$
and a sequence of  positive real numbers $\{\lambda_a: a\in\mathbb{N}\}$
with the following properties:\\
(i) For every $a\in\mathbb{N}$, the function $\psi_a$ satisfies the equation
\begin{equation}\label{5.4}
\frac{4(n+m-1)}{n+m-2}\Delta_{\phi_0}\psi_a-R_{\phi_0}^m \psi_a+\lambda_aw_\infty^{\frac{4}{n+m-2}}\psi_a=0\mbox{ in }M
~~\mbox{ and }~~
\frac{\partial \psi_a}{\partial\nu_{g_0}}=0\mbox{ on }\partial M.
\end{equation}
(ii) For all $a, b\in\mathbb{N}$, we have
\begin{equation}\label{5.5}
\int_M w_\infty^{\frac{4}{n+m-2}}\psi_a\psi_b e^{-\phi_0}dV_{g_0}
=\left\{
   \begin{array}{ll}
     0, & \hbox{if $a\neq b$;} \\
     1, & \hbox{if $a=b$.}
   \end{array}
 \right.
\end{equation}
(iii) The span of $\{\psi_a: a\in\mathbb{N}\}$ is dense in $L^2(M, e^{-\phi_0}dV_{g_0})$.\\
(iv) $\lambda_a\to\infty$ as $a\to\infty$.
\end{proposition}
\begin{proof}
Consider the linear operator $T$
\begin{displaymath}
T: \psi\mapsto w_{\infty}^{-\frac{4}{n+m-2}}\left(\frac{4(n+m-1)}{n+m-2}\Delta_{\phi_0} \psi-R^m_{\phi_0}\psi\right),
\end{displaymath}
where $\psi$ satisfies $\frac{\partial \psi}{\partial\nu_{g_0}}=0$ on $\partial M.$
By integration by parts, we see that this operator $T$ is symmetric with respect to the inner product
\begin{displaymath}
(\psi_1,\psi_2)\mapsto \int_M w_{\infty}^{\frac{4}{n+m-2}} \psi_1\psi_2 e^{-\phi_0}dV_{g_0}
\end{displaymath}
on $L^2(M, e^{-\phi_0}dV_{g_0})$. Hence, the assertion follows from the spectral theorem.
\end{proof}

Let $A$ be a finite subset of $\mathbb{N}$ such that $\lambda_a\leqslant\displaystyle\frac{n+m+2}{n+m-2}r_\infty^m$
for all $a\in A$. We denote by $\Pi$ the projection operator
\begin{equation}\label{5.4A}
\begin{split}
\Pi f&=\sum_{a\not\in A}\left(\int_Mw_\infty^{\frac{4}{n+m-2}}\psi_af e^{-\phi_0}dV_{g_0}\right)w_\infty^{\frac{4}{n+m-2}}\psi_a\\
&=f-\sum_{a\in A}\left(\int_M\psi_af e^{-\phi_0}dV_{g_0}\right)w_\infty^{\frac{4}{n+m-2}}\psi_a.
\end{split}
\end{equation}

In the rest of this section, for simplicity, we denote $W^{1,2}(M, e^{-\phi_0}dV_{g_0})$ and $L^{p}(M, e^{-\phi_0}dV_{g_0})$ by $W^{1,2}(M)$ and $L^{p}(M)$, respectively.
\begin{lemma}\label{lem5.3} For every $1\leqslant p<\infty$, we can find a constant $C$ such that
\begin{equation*}
\begin{split}
\|f\|_{L^p(M)}&\leq C\left\|\frac{4(n+m-1)}{n+m-2}\Delta_{\phi_0}f-R_{\phi_0}^m f+\frac{n+m+2}{n+m-2}r_\infty^mw_\infty^{\frac{4}{n+m-2}}f\right\|_{L^p(M)}\\
&\hspace{4mm}+C\sup_{a\in A}\left|\int_Mw_\infty^{\frac{4}{n+m-2}}\psi_af e^{-\phi_0}dV_{g_0}\right|
\end{split}
\end{equation*}
\end{lemma}
\begin{proof}
Suppose this is not true.
By compactness, we can find a function $f\in L^p(M)$ satisfying
$\|f\|_{L^p(M)}=1$,
\begin{equation}\label{5.6}
\int_M w_\infty^{\frac{4}{n+m-2}}\psi_af e^{-\phi_0}dV_{g_0}=0
\end{equation}
for all $a\in A$ and
\begin{equation*}
\frac{4(n+m-1)}{n+m-2}\Delta_{\phi_0}f-R_{\phi_0}^m f+\frac{n+m+2}{n+m-2}r_\infty^mw_\infty^{\frac{4}{n+m-2}}f=0
\end{equation*}
in the sense of distributions. Hence, if we use $\psi$ as a test function,
then we obtain
$$\left(\lambda_a-\frac{n+m+2}{n+m-2}r_\infty^m\right)\int_Mw_\infty^{\frac{4}{n+m-2}}\psi_a fe^{-\phi_0}dV_{g_0}=0$$
for all $a\in\mathbb{N}$. In particular, we have
$$\int_Mw_\infty^{\frac{4}{n+m-2}}\psi_a fe^{-\phi_0}dV_{g_0}=0$$
for all $a\not\in A$. Combining this with (\ref{5.6}), we can conclude that $f=0$,
which is a contradiction.
\end{proof}

\begin{lemma}\label{lem5.4}
(i) There exists a constant $C$ such that
\begin{equation*}
\begin{split}
\|f\|_{L^{\frac{n+m+2}{n+m-2}}(M)}
&\leqslant C\left\|\Pi\left(\frac{4(n+m-1)}{n+m-2}\Delta_{\phi_0}f-R_{\phi_0}^m f+\frac{n+m+2}{n+m-2}r_\infty^mw_\infty^{\frac{4}{n+m-2}}f\right)\right\|_{L^s(M)}\\
&\hspace{4mm}+C\sup_{a\in A}\left|\int_Mw_\infty^{\frac{4}{n+m-2}}\psi_af e^{-\phi_0}dV_{g_0}\right|
\end{split}
\end{equation*}
where $s=\displaystyle\frac{n(n+m+2)}{n(n+m-2)+2(n+m+2)}.$\\
(ii) There exists a constant such that
\begin{equation*}
\begin{split}
\|f\|_{L^1(M)}
&\leqslant C\left\|\Pi\left(\frac{4(n+m-1)}{n+m-2}\Delta_{\phi_0}f-R_{\phi_0}^m f+\frac{n+m+2}{n+m-2}r_\infty^mw_\infty^{\frac{4}{n+m-2}}f\right)\right\|_{L^1(M)}\\
&\hspace{4mm}+C\sup_{a\in A}\left|\int_Mw_\infty^{\frac{4}{n+m-2}}\psi_af e^{-\phi_0}dV_{g_0}\right|.
\end{split}
\end{equation*}
\end{lemma}
\begin{proof}
If follows from the Sobolev embedding
$W^{2,s}(M)\hookrightarrow L^{\frac{n+m+2}{n+m-2}}(M)$  that
\begin{equation*}
\begin{split}
\|f\|_{L^{\frac{n+m+2}{n+m-2}}(M)}
&\leqslant C\|f\|_{L^s(M)}\\
&\hspace{4mm}
+C\left\| \frac{4(n+m-1)}{n+m-2}\Delta_{\phi_0}f-R_{\phi_0}^m f+\frac{n+m+2}{n+m-2}r_\infty^mw_\infty^{\frac{4}{n+m-2}}f \right\|_{L^s(M)}.
\end{split}
\end{equation*}
This together with Lemma \ref{lem5.3} implies that
\begin{equation*}
\begin{split}
\|f\|_{L^{\frac{n+m+2}{n+m-2}}(M)}
&\leqslant C\sup_{a\in A}\left|\int_Mw_\infty^{\frac{4}{n+m-2}}\psi_af e^{-\phi_0}dV_{g_0}\right|\\
&\hspace{4mm}
+C\left\| \frac{4(n+m-1)}{n+m-2}\Delta_{\phi_0}f-R_{\phi_0}^m f+\frac{n+m+2}{n+m-2}r_\infty^mw_\infty^{\frac{4}{n+m-2}}f \right\|_{L^s(M)}.
\end{split}
\end{equation*}
It follows from the definition of $\Pi$ that
\begin{equation*}
\begin{split}
&\frac{4(n+m-1)}{n+m-2}\Delta_{\phi_0}f-R_{\phi_0}^m f+\frac{n+m+2}{n+m-2}r_\infty^mw_\infty^{\frac{4}{n+m-2}}f\\
&=\Pi\left(\frac{4(n+m-1)}{n+m-2}\Delta_{\phi_0}f-R_{\phi_0}^m f+\frac{n+m+2}{n+m-2}r_\infty^mw_\infty^{\frac{4}{n+m-2}}f\right)\\
&\hspace{4mm}
-\sum_{a\in A}\Big(\lambda_a-\frac{n+m+2}{n+m-2}r_\infty^m\Big)\left(\int_Mw_\infty^{\frac{4}{n+m-2}}\psi_a fe^{-\phi_0}dV_{g_0}\right)w_\infty^{\frac{4}{n+m-2}}\psi_a,
\end{split}
\end{equation*}
which implies that
\begin{equation*}
\begin{split}
&\left\|\frac{4(n+m-1)}{n+m-2}\Delta_{\phi_0}f-R_{\phi_0}^m f+\frac{n+m+2}{n+m-2}r_\infty^mw_\infty^{\frac{4}{n+m-2}}f\right\|_{L^q(M)}\\
&\leqslant \left\|\Pi\left(\frac{4(n+m-1)}{n+m-2}\Delta_{\phi_0}f-R_{\phi_0}^m f+\frac{n+m+2}{n+m-2}r_\infty^mw_\infty^{\frac{4}{n+m-2}}f\right)\right\|_{L^q(M)}\\
&\hspace{4mm}
+C\sup_{a\in A}\left|\int_Mw_\infty^{\frac{4}{n+m-2}}\psi_a fe^{-\phi_0}dV_{g_0}\right|.
\end{split}
\end{equation*}
Now (i) follows from putting these facts together.

Similar to (i), (ii) follows from Lemma \ref{lem5.3} and the definition of $\Pi$.
\end{proof}

\begin{lemma}\label{lem5.5}
There exists a positive real number $\xi$ such that for every vector $z\in\mathbb{R}^A$
with $|z|\leq\xi$, there exists a smooth function $\overline{w}_z$ such that
$\displaystyle\frac{\partial\overline{w}_z}{\partial\nu_{g_0}}=0$ on $\partial M$,
and
\begin{equation*}
\int_Mw_\infty^{\frac{4}{n+m-2}}\psi_a(\overline{w}_z-w_\infty) e^{-\phi_0}dV_{g_0}=z_a
\end{equation*}
for all $a\in A$ and
\begin{equation*}
\Pi\left(\frac{4(n+m-1)}{n+m-2}\Delta_{\phi_0}\overline{w}_z-R_{\phi_0}^m \overline{w}_z+r_\infty^m\overline{w}_z^{\frac{4}{n+m-2}}f\right)=0.
\end{equation*}
Furthermore, the map $z\mapsto\overline{w}_z$ is real analytic.
\end{lemma}
\begin{proof}
This is a consequence of implicit function theorem.
\end{proof}

\begin{lemma}\label{lem5.6}
There exists a real number $0<\gamma<1$ such that
\begin{equation*}
\begin{split}
&E_{(g_0,\phi_0)}(\overline{w}_z)-E_{(g_0,\phi_0)}(w_\infty)\\
&\leq C\sup_{a\in A}\left|\int_M\left(\frac{4(n+m-1)}{n+m-2}\Delta_{\phi_0}\overline{w}_z-R_{\phi_0}^m \overline{w}_z+r_\infty^m\overline{w}_z^{\frac{n+m+2}{n+m-2}}\right)\psi_ae^{-\phi_0}dV_{g_0}\right|^{1+\gamma}
\end{split}
\end{equation*}
if $z$ is sufficiently small.
\end{lemma}
\begin{proof}
Note that the function $z\mapsto E_{(g_0,\phi_0)}(\overline{w}_z)$ is real analytic.
According to the results of Lojasiewicz \cite[equation (2.4)]{Simon83},
there exists a real number $0<\gamma<1$ such that
$$|E_{(g_0,\phi_0)}(\overline{w}_z)-E_{(g_0,\phi_0)}(w_\infty)|\leqslant \sup_{a\in A}\left|\frac{\partial}{\partial a}E(\overline{w}_z)\right|^{1+\gamma}$$
if $z$ is sufficiently small.
For convenience, we define the energy
functional $F_{(g_0,\phi_0)}(w)$ as
\begin{equation*}
F_{(g_0,\phi_0)}(w)
=\frac{\int_MwL_{\phi_0}^m(w)e^{-\phi_0}dV_{g_0}+\int_MwB_{\phi_0}^m(w)e^{-\phi_0}dA_{g_0}}{ \int_Mw^{\frac{2(n+m)}{n+m-2}}e^{-\phi_0}dV_{g_0}}
\end{equation*}
The partial derivatives of the function $z\mapsto E_{(g_0,\phi_0)}(\overline{w}_z)$ are given by
\begin{equation}\label{5.7}
\begin{split}
\frac{\partial}{\partial z_a}E_{(g_0,\phi_0)}(\overline{w}_z)
&=-2\frac{\int_M\left(\frac{4(n+m-1)}{n+m-2}\Delta_{\phi_0}\overline{w}_z-R_{\phi_0}^m\overline{w}_z
+r_\infty^m \overline{w}_z^{\frac{n+m+2}{n+m-2}}\right)\tilde{\psi}_{a,z}e^{-\phi_0}dV_{g_0}}{\left(\int_M\overline{w}_z^{\frac{2(n+m)}{n+m-2}}e^{-\phi_0}dV_{g_0}\right)^{\frac{n+m-2}{n+m}}}\\
&\hspace{4mm}-2(F_{g_0,\phi_0}(\overline{w}_z)-r_\infty^m)
\frac{\int_M\overline{w}_z^{\frac{n+m+2}{n+m-2}}\tilde{\psi}_{a,z}e^{-\phi_0}dV_{g_0}}{\left(\int_M\overline{w}_z^{\frac{2(n+m)}{n+m-2}}e^{-\phi_0}dV_{g_0}\right)^{\frac{n+m-2}{n+m}}}
\end{split}
\end{equation}
where $\tilde{\psi}_{a,z}=\displaystyle\frac{\partial}{\partial z_a}\overline{w}_z$ for $a\in A$.
The function $\tilde{\psi}_{a,z}$ satisfies
\begin{equation*}
\int_Mw_\infty^{\frac{4}{n+m-2}}\tilde{\psi}_{a,z}\psi_b e^{-\phi_0}dV_{g_0}
=\left\{
   \begin{array}{ll}
     0, & \hbox{if $a\neq b$;} \\
     1, & \hbox{if $a=b$}
   \end{array}
 \right.
\end{equation*}
for all $a, b\in A$
and
$$\Pi\left(\frac{4(n+m-1)}{n+m-2}\Delta_{\phi_0}\tilde{\psi}_{a,z}-R_{\phi_0}^m \tilde{\psi}_{a,z}+\frac{n+m+2}{n+m-2}r_\infty^mw_\infty^{\frac{4}{n+m-2}}\tilde{\psi}_{a,z}\right)=0.$$
Using the identity
$$\Pi\left(\frac{4(n+m-1)}{n+m-2}\Delta_{\phi_0}\overline{w}_z-R_{\phi_0}^m \overline{w}_z
+r_\infty^m\overline{w}_z^{\frac{n+m+2}{n+m-2}}\right)=0,$$
we obtain
\begin{equation*}
\begin{split}
&\frac{\partial}{\partial z_a}E_{g_0,\phi_0}(\overline{w}_z)\\
&=-2\frac{\int_M\left(\frac{4(n+m-1)}{n+m-2}\Delta_{\phi_0}\overline{w}_z-R_{\phi_0}^m\overline{w}_z
+r_\infty^m \overline{w}_z^{\frac{n+m+2}{n+m-2}}\right)\psi_a e^{-\phi_0}dV_{g_0}}{\left(\int_M\overline{w}_z^{\frac{2(n+m)}{n+m-2}}e^{-\phi_0}dV_{g_0}\right)^{\frac{n+m-2}{n+m}}}\\
&\hspace{4mm}
+2\sum_{b\in A}\frac{\int_M\overline{w}_z^{\frac{n+m+2}{n+m-2}}\tilde{\psi}_{a,z}e^{-\phi_0}dV_{g_0}}{\left(\int_M\overline{w}_z^{\frac{2(n+m)}{n+m-2}}e^{-\phi_0}dV_{g_0}\right)^{\frac{n+m-2}{n+m}}}\\
&\hspace{8mm}\cdot
\frac{\left(\int_M\left(\frac{4(n+m-1)}{n+m-2}\Delta_{\phi_0}\overline{w}_z-R_{\phi_0}^m\overline{w}_z
+r_\infty^m \overline{w}_z^{\frac{n+m+2}{n+m-2}}\right)\psi_b e^{-\phi_0}dV_{g_0}\right)
\left(\int_M\overline{w}_z^{\frac{4}{n+m-2}}\overline{w}_z\psi_be^{-\phi_0}dV_{g_0}\right)}{\int_M\overline{w}_z^{\frac{2(n+m)}{n+m-2}}e^{-\phi_0}dV_{g_0}}
\end{split}
\end{equation*}
for al $a\in A$. Thus, we obtain
\begin{equation*}
\begin{split}
&\sup_{a\in A}\left|\frac{\partial}{\partial z_a}E_{g_0,\phi_0}(\overline{w}_z)\right|\\
&\leqslant C\sup_{a\in A}\left|\int_M\left(\frac{4(n+m-1)}{n+m-2}\Delta_{\phi_0}\overline{w}_z-R_{\phi_0}^m\overline{w}_z
+r_\infty^m \overline{w}_z^{\frac{n+m+2}{n+m-2}}\right)\psi_a e^{-\phi_0}dV_{g_0}\right|.
\end{split}
\end{equation*}
Combining all these, the assertion follows.
\end{proof}

By Lemma \ref{lem5.5}, the function
\begin{displaymath}
z\to \int_M (w_i-\bar{w}_{z})L^m_{\phi_0}(w_i-\bar{w}_{z})e^{-\phi_0}dV_{g_0}
\end{displaymath}
is analytical and attains the infimum for $|z|\leqslant \xi$. For every $i\in\mathbb{N}$,
we can find $\overline{w}_{z_i}$ such that $|z_i|\leq \xi$ and
$$\int_M(w_i-\overline{w}_{z_i})L_{\phi_0}^m(w_i-\overline{w}_{z_i})e^{-\phi_0}dV_{g_0}
\leq \int_M(w_i-\overline{w}_{z})L_{\phi_0}^m(w_i-\overline{w}_{z})e^{-\phi_0}dV_{g_0}$$
for all $|z|\leq \xi$.

Notice that by Lemma \ref{lem5.5}, we have $\overline{w}_{0}=w_{\infty}$. Combining this with the definition of $\overline{w}_{z_i}$, we have
$$\int_M(w_i-\overline{w}_{z_i})L_{\phi_0}^m(w_i-\overline{w}_{z_i})e^{-\phi_0}dV_{g_0}
\leq \int_M(w_i-w_\infty)L_{\phi_0}^m(w_i-w_\infty)e^{-\phi_0}dV_{g_0}.$$

By the compactness result in Proposition \ref{prop5.2}, the
expression on the right hand side tends to zero as $i\to\infty$,
i.e. we have as $i\to\infty$,
\begin{equation}\label{5.8}
\|w_i-\overline{w}_{z_i}\|_{W^{1,2}(M)}\to 0~~\mbox{ and }~~\|\overline{w}_{z_i}-w_\infty\|_{W^{1,2}(M)}\to 0.
\end{equation}

\begin{lemma}\label{lem5.7}
The difference $w_i-\overline{w}_{z_i}$ satisfies
$$\|w_i-\overline{w}_{z_i}\|_{L^{\frac{n+m+2}{n+m-2}}(M)}
\leqslant C\left\|w_i^{\frac{n+m+2}{n+m-2}}(R_{\phi_i}^m-r_\infty^m)\right\|_{L^{\frac{2(n+m)}{n+m-2}}(M)}
+o(1)$$
if $i$ is sufficiently large.
\end{lemma}
\begin{proof}
For simplicity, we denote $w_i-\bar{w}_{z_i}$ by $u_i$. Using the identities
\begin{equation}\label{5.10}
\frac{4(n+m-1)}{n+m-2}\Delta_{\phi_0} w_i-R^m_{\phi_0}w_i+r^m_{\infty}w_i^{\frac{n+m+2}{n+m-2}}=-w_i^{\frac{n+m+2}{n+m-2}}(R^m_{\phi_i}-r^m_{\infty})
\end{equation}
and
\begin{displaymath}
\Pi\left( \frac{4(n+m-1)}{n+m-2}\Delta_{\phi_0} \bar{w}_z-R^m_{\phi_0}\bar{w}_z+r^m_{\infty}\bar{w}_z^{\frac{n+m+2}{n+m-2}}\right)=0,
\end{displaymath}
we obtain
\begin{equation}
\begin{aligned}
\Pi&\left(\frac{4(n+m-1)}{n+m-2}\Delta_{\phi_0}u_i-R^m_{\phi_0}u_i+\frac{n+m+2}{n+m-2}r^m_{\infty}w_{\infty}^{\frac{4}{n+m-2}}u_i\right)\\
=& \Pi\left(-w_i^{\frac{n+m+2}{n+m-2}}(R^m_{\phi_i}-r^m_{\infty})-\frac{n+m+2}{n+m-2}r^m_{\infty}\left(\bar{w}_{z_i}^{\frac{4}{n+m-2}}-w_{\infty}^{\frac{4}{n+m-2}}\right)u_i \notag\right.
\\
\phantom{=\;\;}
& \left. +r^m_{\infty}\left(w_i^{\frac{n+m+2}{n+m-2}}-\bar{w}_{z_i}^{\frac{n+m+2}{n+m-2}}+ \frac{n+m+2}{n+m-2}\bar{w}_{z_i}^{\frac{4}{n+m-2}}u_i\right)\right)
\end{aligned}
\end{equation}
It follows from Lemma \ref{lem5.4}(i) that
\begin{displaymath}
\begin{split}
    &\| u_i\|_{L^{\frac{n+m+2}{n+m-2}}(M)}\leqslant  C\sup_{a\in A}\left|\int_M w_{\infty}^{\frac{4}{n+m-2}}\psi_a u_i
    e^{-\phi_0}{\rm dvol}_{g_0}\right|\\
    & +C\left\| \Pi\left(\frac{4(n+m-1)}{n+m-2}\Delta_{ \phi_0}u_i -R^m_{\phi_0}u_i+\frac{n+m+2}{n+m-2}r^m_{\infty}w_{\infty}^{\frac{4}{n+m-2}}u_i\right)\right\|_{L^{s}(M)}.
\end{split}
\end{displaymath}
We conclude that
\begin{displaymath}
\begin{split}
    \|u_i\|&_{L^{\frac{n+m+2}{n+m-2}}(M)}\leqslant  C\left\| \bar{w}_{z_i}^{\frac{n+m+2}{n+m-2}}-w_i^{\frac{n+m+2}{n+m-2}}+ \frac{n+m+2}{n+m-2}\bar{w}_{z_i}^{\frac{4}{n+m-2}}u_i \right\|_{L^{s}(M)}\\
    & +C\left\| w_i^{\frac{n+m+2}{n+m-2}}(R^m_{\phi_i}-r^m_{\infty})\right\|_{L^{s}(M)}+C\left\| (\bar{w}_{z_i}^{\frac{4}{n+m-2}}-w_{\infty}^{\frac{4}{n+m-2}})u_i \right\|_{L^{s}(M)}\\
    & +C\sup_{a\in A}\left|\int_M w_{\infty}^{\frac{4}{n+m-2}}\psi_a u_i
    e^{-\phi_0}dV_{g_0}\right|.
\end{split}
\end{displaymath}
According to (\ref{5.8}), without loss of generality, we can assume that
\begin{equation}\label{pointcong}
    w_i\to w_{\infty} ~~\mbox{ and }~~ \bar{w}_{z_i}\to w_{\infty} \quad ~~\mbox{ a.e. in }~~ M.
\end{equation}
Combining (\ref{pointcong}) with Lebesgue's dominated convergence theorem yields
\begin{equation*}
    \left\| (\bar{w}_{z_i}^{\frac{4}{n+m-2}}-w_{\infty}^{\frac{4}{n+m-2}})u_i \right\|_{L^{s}(M)}=o(1).
\end{equation*}

By (\ref{pointcong}), we have the pointwise estimate
\begin{equation*}
\left|\bar{w}_{z_i}^{\frac{n+m+2}{n+m-2}}-w_i^{\frac{n+m+2}{n+m-2}}+ \frac{n+m+2}{n+m-2}\bar{w}_{z_i}^{\frac{4}{n+m-2}}u_i \right|
    \leqslant C \bar{w}_{z_i}^{\frac{n+m+2}{n+m-2}-2}|u_i|^{2}
\end{equation*}
if $i$ is sufficiently large. Since $w_{\infty}$ is a positive smooth function in $M$, we obtain
\begin{equation*}
    \left\| \bar{w}_{z_i}^{\frac{n+m+2}{n+m-2}}-w_i^{\frac{n+m+2}{n+m-2}}+ \frac{n+m+2}{n+m-2}\bar{w}_{z_i}^{\frac{4}{n+m-2}}u_i \right\|_{L^{s}(M)} \leqslant C \left\| u_i \right\|^{2}_{L^{2s}(M)}.
\end{equation*}
When $n\geqslant 3$ and $m>0$, $2s=\frac{2n(n+m+2)}{n(n+m-2)+2(n+m+2)}<\frac{2(n+m)}{n+m-2}$. By H\"older's inequality, we conclude that
\begin{equation*}
    \left\| \bar{w}_{z_i}^{\frac{n+m+2}{n+m-2}}-w_i^{\frac{n+m+2}{n+m-2}}+ \frac{n+m+2}{n+m-2}\bar{w}_{z_i}^{\frac{4}{n+m-2}}u_i \right\|_{L^{s}(M)}=o(1).
\end{equation*}
Moreover, since the set $A$ is a finite, we have
\begin{equation*}
    \sup_{a\in A}\left|\int_M w_{\infty}^{\frac{4}{n+m-2}}\psi_a u_i
    e^{-\phi_0}dV_{g_0}\right|\leqslant C\left\|u_i\right\|_{L^1(M)}=o(1).
\end{equation*}
Putting these facts together, the assertion follows.
\end{proof}

\begin{lemma}\label{lem5.8}
The difference $w_i-\overline{w}_{z_i}$ satisfies
$$\|w_i-\overline{w}_{z_i}\|_{L^1(M)}
\leqslant C\left\|w_i^{\frac{n+m+2}{n+m-2}}(R_{\phi_i}^m-r_\infty^m)\right\|_{L^{\frac{2(n+m)}{n+m-2}}(M)}
+o(1)$$
if $i$ is sufficiently large.
\end{lemma}
\begin{proof}
The proof is almost the same as that of Lemma \ref{lem5.7},
except we use Lemma \ref{lem5.4}(ii) instead of Lemma \ref{lem5.4}(i).
We omit the proof and leave it to the readers.
\end{proof}

\begin{lemma}\label{lem5.9} There holds
\begin{equation*}
\begin{split}
&\sup_{a\in A}\left|\int_M\left(\frac{4(n+m-1)}{n+m-2}\Delta_{\phi_0}\overline{w}_{z_i}-R_{\phi_0}^m \overline{w}_{z_i}
+r_\infty^m\overline{w}_{z_i}^{\frac{n+m+2}{n+m-2}}\right)\psi_a e^{-\phi_0}dV_{g_0}\right|\\
&\leqslant C\left(\int_M w(t_i)^{\frac{2(n+m)}{n+m-2}}|R_{\phi(t_i)}^m-r_\infty^m|^{\frac{2(n+m)}{n+m+2}}e^{-\phi_0}dV_{g_0}\right)^{\frac{n+m+2}{2(n+m)}}+o(1)
\end{split}
 \end{equation*}
if $i$ is sufficiently large.
\end{lemma}
\begin{proof}
Integration by parts yields
\begin{equation*}
\begin{split}
&\int_M\left(\frac{4(n+m-1)}{n+m-2}\Delta_{\phi_0}\overline{w}_{z_i}-R_{\phi_0}^m \overline{w}_{z_i}
+r_\infty^m\overline{w}_{z_i}^{\frac{n+m+2}{n+m-2}}\right)\psi_a e^{-\phi_0}dV_{g_0}\\
&=\int_M\left(\frac{4(n+m-1)}{n+m-2}\Delta_{\phi_0}w_i-R_{\phi_0}^m w_i
+r_\infty^mw_i^{\frac{n+m+2}{n+m-2}}\right)\psi_a e^{-\phi_0}dV_{g_0}\\
&\hspace{4mm}+\lambda_a\int_Mw_\infty^{\frac{4}{n+m-2}}(w_i-\overline{w}_{z_i})e^{-\phi_0}dV_{g_0}
-r_\infty^m\int_M\left(w_i^{\frac{n+m+2}{n+m-2}}-\overline{w}_{z_i}^{\frac{n+m+2}{n+m-2}}\right)\psi_a e^{-\phi_0}dV_{g_0},
\end{split}
\end{equation*}
where we have used the fact that
$\displaystyle\frac{\partial \overline{w}_{z_i}}{\partial\nu_{g_0}}=\frac{\partial w_i}{\partial\nu_{g_0}}=\frac{\partial \psi_a}{\partial\nu_{g_0}}=0$
on $\partial M$.
Combining this with (\ref{5.10}) yields
\begin{equation*}
\begin{split}
&\int_M\left(\frac{4(n+m-1)}{n+m-2}\Delta_{\phi_0}\overline{w}_{z_i}-R_{\phi_0}^m \overline{w}_{z_i}
+r_\infty^m\overline{w}_{z_i}^{\frac{n+m+2}{n+m-2}}\right)\psi_a e^{-\phi_0}dV_{g_0}\\
&=-\int_Mw_i^{\frac{n+m+2}{n+m-2}}(R_{\phi_i}^m-r_\infty^m)\psi_a e^{-\phi_0}dV_{g_0}\\
&\hspace{4mm}+\lambda_a\int_Mw_\infty^{\frac{4}{n+m-2}}(w_i-\overline{w}_{z_i})e^{-\phi_0}dV_{g_0}
-r_\infty^m\int_M\left(w_i^{\frac{n+m+2}{n+m-2}}-\overline{w}_{z_i}^{\frac{n+m+2}{n+m-2}}\right)\psi_a e^{-\phi_0}dV_{g_0}.
\end{split}
\end{equation*}
Using the pointwise estimate
$$\left|w_i^{\frac{n+m+2}{n+m-2}}-\overline{w}_{z_i}^{\frac{n+m+2}{n+m-2}}\right|
\leq C\overline{w}_{z_i}^{\frac{4}{n+m-2}}|w_i-\overline{w}_{z_i}|+C|w_i-\overline{w}_{z_i}|^{\frac{n+m+2}{n+m-2}},$$
we can then deduce that
\begin{equation*}
\begin{split}
&\sup_{a\in A}\left|\int_M\left(\frac{4(n+m-1)}{n+m-2}\Delta_{\phi_0}\overline{w}_{z_i}-R_{\phi_0}^m \overline{w}_{z_i}
+r_\infty^m\overline{w}_{z_i}^{\frac{n+m+2}{n+m-2}}\right)\psi_a e^{-\phi_0}dV_{g_0}\right|\\
&\leq C\left(\int_Mw(t_i)^{\frac{2(n+m)}{n+m-2}}|R_{\phi(t_i)}^m-r_\infty^m|^{\frac{2(n+m)}{n+m-2}} e^{-\phi_0}dV_{g_0}\right)^{\frac{n+m+2}{2(n+m)}}\\
&\hspace{4mm}+C\|w_i-\overline{w}_{z_i}\|_{L^1(M)}+
C\|w_i-\overline{w}_{z_i}\|_{L^{\frac{n+m+2}{n+m-2}}(M)}^{\frac{n+m+2}{n+m-2}}.
\end{split}
\end{equation*}
Now the assertion follows from
combining this with Lemmas \ref{lem5.7} and \ref{lem5.8}.
\end{proof}

Combining Lemma \ref{lem5.6} and Lemma \ref{lem5.9},
we immediately have the following:

\begin{proposition}\label{prop5.11} We have the following estimate
\begin{equation*}
\begin{split}
&E_{(g_0,\phi_0)}(\overline{w}_z)-E_{(g_0,\phi_0)}(w_\infty)\\
&\leqslant C\left(\int_M w(t_i)^{\frac{2(n+m)}{n+m-2}}|R_{\phi(t_i)}^m-r_\infty^m|^{\frac{2(n+m)}{n+m+2}}e^{-\phi_0}dV_{g_0}\right)^{\frac{n+m+2}{2(n+m)}(1+\gamma)}+o(1)
\end{split}
\end{equation*}
if $i$ is sufficiently large.
\end{proposition}

We are now ready to prove Proposition \ref{prop4.3}.

\begin{proof}[Proof of Proposition \ref{prop4.3}]
It follows from the definition of $r_\phi^m$
and the assumption (\ref{vol}) that
$$r_{\phi(t_i)}^m-r_\infty^m=E_{(g_0,\phi_0)}(w_i)-E_{(g_0,\phi_0)}(w_\infty).$$
Moreover, it follows from (\ref{5.8}) that
$$E_{(g_0,\phi_0)}(w_i)=E_{(g_0,\phi_0)}(\overline{w}_{z_i})+o(1).$$
Now the assertion combining all these with
Proposition \ref{prop5.11}.
\end{proof}

\section{Nonpositive cases}

In this section, we deal with the remaining cases; i.e. $Y_{n,m}[(g_0, \phi_0)]\leqslant 0$.

\subsection{Negative case}
As discussed in Lemma \ref{threecases}, we can choose an initial metric measure space $(M, g_0, e^{-\phi_0}dV_{g_0},e^{-\phi_0}dA_{g_0}, m)$ such that
\begin{equation}
R^m_{\phi_0}<0 ~~\mbox{ in }M~~\mbox{ and }~~H_{g}=0\mbox{ on }\partial M.
\end{equation}
Let $w(t)$ be the solution of (\ref{wevolution}) on a maximal time interval $[0,T^*)$. Applying the maximal principle to (\ref{wevolution}) derives
\begin{equation}\label{5.2A}
    \frac{d }{d t}w_{\min}^{N}(t)\geqslant \frac{n+m+2}{4}\left( \min |R^m_{\phi_0}|w_{\min}(t)+r^m_{\phi}w_{\min}^{N}(t)\right),
\end{equation}
where $w_{\min}(t)=\displaystyle\min_{M}w(t)$ and $N=\frac{n+m+2}{n+m-2}$. By  (\ref{normalize2}), we have
\begin{equation}\label{rY}
    r^m_{\phi(t)}\geqslant Y_{n,m}[(g_0, \phi_0)].
\end{equation}
Note that $Y_{n,m}[(g_0, \phi_0)]$ is finite by H\"older's inequality. Hence, integrating (\ref{5.2}) yields
\begin{equation}\label{neglow}
    w_{\min}^{N-1}(t)\geqslant C\cdot\min \left\{w_{\min}^{N-1}(0), \frac{\min |R^m_{\phi_0}|}{|Y_{n,m}[(g_0, \phi_0)]|}\right\}
\end{equation}
for some uniform constant $C$.
On the other hand, applying the maximal principle to (\ref{wevolution}) also gives
\begin{equation}\label{5.5A}
    \frac{d }{d t}w_{\max}^{N}(t)\leqslant \frac{n+m+2}{4}\left( -\left(\min R^m_{\phi_0}\right)w_{\max}(t)+r^m_{\phi}w_{\max}^{N}(t)\right),
\end{equation}
where $w_{\max}(t)=\displaystyle\max_{M}w(t)$. According to (\ref{eq:r}), we obtain that
\begin{equation}\label{negsup}
    w_{\max}^{N}(t)\leqslant \left(w_{\max}^{N}(0)+1\right)e^{c(|\min R^m_{\phi_0}|+r^m_{\phi(0)})t},
\end{equation}
for some positive constant $c$. It follows from (\ref{neglow}) and (\ref{negsup})   that $w(t)$ will not blow up in finite time; i.e. $T^*=\infty$.

We claim that $r^m_{\phi(t)}$ will eventually become negative, even if $r^m_{\phi(0)}$ may not be so. If $r^m_{\phi(t)}$ is always nonnegative for $t\geqslant 0$, (\ref{5.2}) would imply
\begin{equation} 
    \frac{d }{d t}w_{\min}^{N}(t)\geqslant \frac{n+m+2}{4} \min |R^m_{\phi^0}|w_{\min}(t).
\end{equation}
Hence $w_{\min}(t)$ approaches to infinity as $t\to \infty$, which contradicts (\ref{normalize1}). Without loss of generality, we may assume $r^m_{\phi(0)}<0$. 
By (\ref{5.5A}), we have
\begin{equation*}
    w_{\max}^{N-1}(t)\leqslant C\cdot \max\left\{w_{\max}^{N-1}(0), \frac{\max |R^m_{\phi^0}|}{|r^m_{\phi(0)}|}\right\}
\end{equation*}
for some uniform positive constant $C$. This together with (\ref{neglow}) implies that $w(t)$ is uniformly bounded from above and away from zero.

Moreover, it follows from (\ref{Revolution}) that
\begin{equation*}
\frac{d (R^m_{\phi})_{\min}}{d t}\geqslant (R^m_{\phi})_{\min}((R^m_{\phi})_{\min}-r^m_{\phi})\geqslant r^m_{\phi}((R^m_{\phi})_{\min}-r^m_{\phi}),
\end{equation*}
where $(R^m_{\phi})_{\min}=\displaystyle\min_{M}R^m_{\phi}(t)$.
Combining this with (\ref{rY}), we can obtain a uniform lower bound on $ R^m_{\phi}(t)$; i.e. for all $t\geqslant 0$
\begin{equation}
   R^m_{\phi}(t)\geqslant r^m_{\phi}(t)-Ce^{r^m_{\phi(0)}t}\geqslant Y_{n,m}[(g_0, \phi_0)]-C.
\end{equation}
Similar to Proposition \ref{prop3.4}, the maximum principle also implies
\begin{equation*}
     \sup_{M} R^m_{\phi}(t)\leqslant \max \left\{ \sup_{M}R^m_{\phi}(0), 0 \right\}
\end{equation*}
Therefore, we can generalize Lemma \ref{lem3.7} and Proposition \ref{Holder} to the negative case. In view of the argument at the end of Section 2, we 
can derive uniform estimates for all higher order derivatives of $w(t)$, $t\geqslant 0$.

\subsection{Zero case}
Here, we treat the zero case. 
As discussed in Lemma \ref{threecases}, 
we can fix a background metric measure space $(M, g_0, e^{-\phi_0}dV_{g_0}, e^{-\phi_0}dA_{g_0},m)$ such that $R^m_{\phi_0}\equiv 0$. Note that by \cite[Proposition 3.5]{Case15}, $r^m_{\phi(t)}$ can never be negative. Since the function $t\mapsto r^m_{\phi}(t)$ is nonincreasing, $r^m_{\phi(0)}=0$ implies $r^m_{\phi(t)}\equiv 0$. Thus the solution of (\ref{flow}) is constant in time.

We next assume that $r^m_{\phi(0)}>0$. Observe that
\begin{equation*}
    \frac{w_{\min}^{N}(t)}{w_{\min}^{N}(0)}\geqslant c \int_{0}^t r^m_{\phi(t)} dt ~~\mbox{ and }~~ \frac{w_{\max}^{N}(t)}{w_{\max}^{N}(0)}\leqslant c \int_{0}^t r^m_{\phi(t)} dt
\end{equation*}
for some positive constant $c$, which are the consequences of (\ref{5.2A}) and (\ref{5.5A}). Hence we obtain the Harnack following inequality
\begin{equation*}
     \frac{w_{\min}^{N}(t)}{w_{\min}^{N}(0)} \geqslant \frac{w_{\max}^{N}(t)}{w_{\max}^{N}(0)}.
\end{equation*}
It follows that $w(t)$ exists for all time.
By the same argument as in Subsection 5.1, we can derive the smooth convergence.

\section{Acknowledgement}

The first author was supported
by Basic Science Research Program through the National Research Foundation of Korea (NRF) funded by the Ministry of Education, Science and Technology (Grant No. 2020R1A6A1A03047877 and 2019R1F1A1041021), and by Korea Institute for Advanced Study (KIAS) grant
funded by the Korea government (MSIP). The second author was supported by a KIAS Individual Grant (MG070701) at Korea Institute for Advanced Study.

\bibliographystyle{amsplain}

\end{document}